\documentclass[onecolumn,final,a4page]{siamltex}

\usepackage{amsmath}
\usepackage{amssymb}
\usepackage{graphicx}
\usepackage[abs]{overpic}
\usepackage{subfigure}
\usepackage{color}
\usepackage{epsfig}
\usepackage{wrapfig}
\usepackage{comment}
\usepackage[makeroom]{cancel}
\usepackage{bm}
\usepackage{rotating}
\usepackage[numbers,sort]{natbib}
\usepackage{textpos}

\newcommand{\email}[1]{\emph{#1}}

\newtheorem{remark}{Remark}

\newtheorem{assum}{Assumption}
\newtheorem{thm}{Theorem}
\newtheorem{lem}{Lemma}

\def\Div{\operatorname{div}}

\def\d{\partial}

\def\D{\nabla}
\def\R{\mathbb{R}}

\def\Bq{\mathbf{q}}
\def\Bw{\mathbf{w}}

\def\Bu{\mathbf{u}}

\def\Bf{\mathbf{f}}
\def\Bg{\mathbf{g}}

\def\Bx{\mathbf{x}}

\def\Bv{\mathbf{v}}

\def\Bn{\mathbf{n}}
\def\Bm{\mathbf{m}}             

\def\BW{\mathbf{W}}

\def\cE{\mathcal{E}}

\def\cM{\mathcal{M}}

\def\mI{\mathbb{I}}

\def\mD{\mathbb{D}}

\def\mH{\mathbb{H}}

\newcommand{\hsk}{\hspace{-0.25mm}}
\newcommand{\TNORM} [1]{|\hsk|\hsk|#1|\hsk|\hsk|}
\newcommand{\TNORME}[1]{|\hsk|\hsk|#1|\hsk|\hsk|_{E}}

\newcommand{\bSnorm} [1]{\norm{#1}_{1,h}}
\newcommand{\bSnormE}[1]{\norm{#1}_{1,E}}

\title{The non-conforming Virtual Element Method for\\the Stokes
    equations
    }

\author{
  Andrea Cangiani\footnotemark[1]\,, 
  Vitaliy Gyrya\footnotemark[2]\,,\!\! 
  \and 
  Gianmarco Manzini\footnotemark[2]\ \footnotemark[3]
}
  
\begin{document}

\maketitle

\pagestyle{myheadings}
\thispagestyle{plain}
\markboth{A.~Cangiani, V.~Gyrya, and G.~Manzini}
{The Non-conforming Virtual Element Method for the Stokes equations}

\renewcommand{\thefootnote}{\fnsymbol{footnote}}

\footnotetext[1]{
  Department of Mathematics, University of Leicester,
  University Road - Leicester LE1 7RH, United Kingdom 
  \email{andrea.cangiani@univ-leicester.uk}
}
\footnotetext[2]{
Los Alamos National Laboratory, Theoretical Division,
    Group T-5, MS B284, Los Alamos, NM-87545, USA
    \email{\{gyrya,gmanzini\}@lanl.gov}
}

\footnotetext[3]{
  Istituto di Matematica Applicata e Tecnologie
  Informatiche (IMATI) -- CNR, via Ferrata 1, I -- 27100 Pavia, Italy
}

\renewcommand{\thefootnote}{\arabic{footnote}}

\begin{abstract}
  We present the non-conforming Virtual Element Method (VEM) for the
  numerical approximation of velocity and pressure in the steady 
  Stokes problem.
  The pressure is approximated using discontinuous piecewise
  polynomials, while each component of the velocity is approximated
  using the nonconforming virtual element space.
  On each mesh element the local virtual space contains the 
  space of polynomials of up to a given degree, plus suitable non-polynomial 
  functions. The virtual element functions are implicitly
  defined as the solution of local Poisson problems with 
  polynomial Neumann boundary conditions.
  As typical in VEM approaches,
  the explicit evaluation of the non-polynomial functions
  is not required.
  This approach makes it possible to construct nonconforming (virtual)
  spaces for any polynomial degree regardless 
  of the parity, for two-and three-dimensional problems, and for 
  meshes with very general polygonal and polyhedral elements.
  We show that the non-conforming VEM is inf-sup stable
  and establish optimal a priori error estimates for the velocity 
  and pressure approximations.
  Numerical examples confirm the convergence 
  analysis and the effectiveness of the method in providing high-order 
  accurate approximations.
\end{abstract}




\raggedbottom

\newcommand{\dx}{d\mathbf{x}}
\newcommand{\ds}{ds}

\newcommand{\Edges}{\cE_h}
\newcommand{\Th}{\mathcal{T}_h}

\newcommand{\uh} {u_h}
\newcommand{\vh} {v_h}
\newcommand{\wh} {w_h}
\newcommand{\ph} {p_h}
\newcommand{\qh} {q_h}
\newcommand{\qI} {q^{\INTP}}
\newcommand{\pI} {p^{\INTP}}

\newcommand{\Buh}{\Bu_h}
\newcommand{\Bvh}{\Bv_h}
\newcommand{\Bwh}{\Bw_h}
\newcommand{\Bfh}{\Bf_h}
\newcommand{\BuI}  {\Bu^{\INTP}}
\newcommand{\BvI}  {\Bv^{\INTP}}
\newcommand{\Bvq}  {\Bv_q}
\newcommand{\Bvtq} {\hat{\Bv}_q}
\newcommand{\BvtqI}{\hat{\Bv}_q^{\INTP}}
\newcommand{\Bgh}{\Bg_h}
\newcommand{\Bvhb}{\overline{\Bv}_h}

\newcommand{\as}{a}
\newcommand{\bs}{b}

\newcommand{\aE}{a^E}
\newcommand{\bE}{b^E}

\newcommand{\ah}{a_h}
\newcommand{\bh}{b_h}

\newcommand{\ahE}{a_h^E}
\newcommand{\bhE}{b_h^E}

\newcommand{\sh} {S_h}
\newcommand{\shE}{S_h^E}

\newcommand{\PE}[1]{\mathbb{P}_{#1}(E)}
\newcommand{\PS}[1]{\mathbb{P}_{#1}}

\newcommand{\Prod}{\Pi}
\newcommand{\PiE}[1]{\Pi^{E}_{#1}}
\newcommand{\Pis}[1]{\Pi^{s}_{#1}}

\newcommand{\INTP} {\mathtt{I}}

\newcommand{\VhE}{V_{E}^k}
\newcommand{\Whk}{\BW_{h}^k}
\newcommand{\Vhk}{V_{h}^k}
\newcommand{\Phk}{\Phi_{h}^{k-1}}
\newcommand{\Vhkz}{V_{h,0}^k}
\newcommand{\Vhkg}{V_{h,g}^k}
\newcommand{\BPk}{\mathbf{P}^{k}_{h}}
\newcommand{\HDIV}[1]{H(\mbox{div},#1)}
\newcommand{\VVhk}{{\mathbf{V}}_{h}^k}    
\newcommand{\VVhkz}{{\mathbf{V}}_{h,0}^k} 
\newcommand{\VVhkg}{{\mathbf{V}}_{h,g}^k} 

\newcommand{\reg}{m}

\newcommand{\Bma}{\Bm_{\alpha}}
\newcommand{\ma}{m_{\alpha}}
\newcommand{\mb}{m_{\beta}}

\newcommand{\Abs}  [1]{{\Big\lvert#1\Big\rvert}} 
\newcommand{\abss} [1]{{\big\lvert#1\big\rvert}} 
\newcommand{\abs}  [1]{{\left\lvert#1\right\rvert}} 
\newcommand{\norm} [1]{{\left\lVert#1\right\rVert}} 
\newcommand{\tnorm}[1]{\lvert\hspace{-0.3mm}\lvert\hspace{-0.3mm}\lvert#1\rvert\hspace{-0.3mm}\rvert\hspace{-0.3mm}\rvert}      
\newcommand{\jump} [1]{{\left[\!\left[#1\right]\!\right]}}
\newcommand{\jumpt}[1]{\widetilde{{\left[\!\left[#1\right]\!\right]}}}
\newcommand{\Id}{\operatorname{I}} 
\newcommand{\Real}{\mathbb{R}}

\newcommand{\MAG}[1]{{\color{magenta}#1}}
\newcommand{\RED}[1]{{\color{red}#1}}
\newcommand{\EOD}{\end{document}}

\def\trait #1 #2 #3 {\vrule width #1pt height #2pt depth #3pt}
\def\fin{\hfill
        \trait .3 5 0
        \trait 5 .3 0
        \kern-5pt
        \trait 5 5 -4.7
        \trait 0.3 5 0
\medskip}
\newcommand{\ENDPROOF}{\fin}
\newcommand{\BEGINPROOF}{~\\\emph{Proof}.~~}


\section{Introduction}

We are concerned with the development of the non-conforming virtual
element method (VEM) for the Stokes problem in the unknown fields
$\Bu$ and $p$ satisfying
\begin{align}
  -\Delta\Bu + \D p &= \Bf\phantom{0\Bg}   \mbox{in~}\Omega,\label{eq:Stokes:PDE:a} \\
  \Div\Bu           &= 0  \phantom{\Bf\Bg} \mbox{in~}\Omega,\label{eq:Stokes:PDE:b} \\
  \Bu               &= \Bg\phantom{0\Bf}   \mbox{on~}\Gamma,\label{eq:Stokes:PDE:c}
\end{align}
where $\Omega$ is a polygonal or polyhedral domain in $\R^d$, $d=2,3$
with boundary $\Gamma$.
We will refer to $\Bu$ and $p$ as \emph{velocity} and \emph{pressure},
respectively.

Historically, the first non-conforming finite element space dates back
to the work of Crouzeix and Raviart in~\cite{Crouzeix-Raviart:1973}.
Their method provides a low-order accurate approximation of the
velocity field of the Stokes equations on triangular meshes based on
linear polynomials.
Later on, higher-order accurate methods were proposed by Fortin and
Soulie~\cite{Fortin-Soulie:1983}, and Crouzeix and Falk
\cite{Crouzeix-Falk:1989}, respectively, by using finite element
spaces based on polynomials of degree $k=2$ and $3$ on triangles.
The functions in these finite element spaces are continuous on a
discrete set of points located at the internal mesh edges.
These points are the roots of the one-dimensional $k^{th}$-order
Legendre polynomials defined over the edges and can be used as the
nodes of the $k^{th}$ Gauss-Legendre quadrature rule.
This minimal continuity requirement ensures the optimal convergence
rate; see, for instance,~\cite{Crouzeix-Raviart:1973}.
The construction of the non-conforming elements for the Stokes problem
has been recently generalized on triangles in
\cite{Stoyan-Baran:2006,Baran-Stoyan:2007} to consider polynomials of
any degree $k>3$, thus resulting in the so-called family of
Gauss-Legendre non-conforming methods.
Robust a posteriori estimates for such schemes can be found in
\cite{Ainsworth:2005}.
A major drawback of the non-conforming Gauss-Legendre elements is that
the space construction for even $k$ differs from that of odd $k$.
This feature also affects the classical low-order cases for $k\leq 3$,
i.e., the formulation of the non-conforming spaces for $k=1$ and $k=3$
in \cite{Crouzeix-Raviart:1973,Crouzeix-Falk:1989} is not the same as
for $k=2$ in~\cite{Fortin-Soulie:1983}.

The generalization of the non-conforming formulation to elements other
than triangles in two and three dimensions is quite a hard task due to
the difficulty of the construction of the shape functions for such
elements.
For example, successful attempts in this direction are found for
quadrilaterals, tetrahedra and hexahedra
in~\cite{Rannacher-Turek:1992,cai1999stable,Matthies-Tobiska:2005,Matthies:2007}.
Instead, the construction of the non-conforming virtual element space
for the Stokes equations that we present in this paper is
straightforward for any polynomial degree regardless of its parity,
and is the same for elements with very general geometric shape in two
and three space dimensions.

The VEM was first introduced as a $C^0$-conforming formulation for the
Poisson equation with constant coefficients
in~\cite{BeiraodaVeiga-Brezzi-Cangiani-Manzini-Marini-Russo:2013}.
The non-conforming formulation for the same problem was developed
later
in~\cite{Ayuso-Lipnikov-Manzini:2016,Cangiani-Manzini-Sutton:2016}.
In both formulations, the trial and test functions are defined
implicitly on each mesh element as the solution of a boundary value
problem and never explicitly constructed in practice, hence the name
``\emph{virtual}''.
In view of obtaining a \emph{computable} and \emph{accurate} virtual
formulation, two essential ingredients are sought for the virtual
element space: $(i)$ it must contain a space of polynomials up to a
given degree; $(ii)$ orthogonal $L^2$ and $H^1$ projections of virtual
functions onto the polynomial sub-space must be computable just using
the degrees of freedom.
Properties $(i)$-$(ii)$ make possible to avoid the explicit
construction of the shape functions and allows us to formulate and
implement the method on very general polygonal and polyhedral meshes.
These features are inherited from the Mimetic Finite Difference (MFD)
method~\cite{Lipnikov-Manzini-Shashkov:2014,BeiraodaVeiga-Lipnikov-Manzini:2014},
which can be seen as a precursor.
A conforming low order MFD method for the steady Stokes problem on
general polygonal and polyhedral meshes is found
in~\cite{BeiraodaVeiga-Gyrya-Lipnikov-Manzini:2009,BeiraodaVeiga-Lipnikov-Manzini:2010},
and a higher order MFD method equivalent to the non-conforming VEM
in~\cite{Ayuso-Lipnikov-Manzini:2016} is found
in~\cite{Lipnikov-Manzini:2014:JCP}.

The objective of this paper is to develop the non-conforming VEM for
the weak form of~\eqref{eq:Stokes:PDE:a}-\eqref{eq:Stokes:PDE:c}
(see~\eqref{eq:Stokes:PDE:weak:a}-\eqref{eq:Stokes:PDE:weak:b} in the
next section) that is suitable for very general mesh partitioning of
$\Omega$ in polygons and polyhedra.
As is standard in the finite element setting, the VEM approximation of
\eqref{eq:Stokes:PDE:a}-\eqref{eq:Stokes:PDE:c} proposed in this work
is based on the construction of a pair of finite element spaces
satisfying the inf-sup condition, see~\cite{Boffi-Brezzi-Fortin:2013}.
The major features of this VEM are:
$(i)$ each component of the velocity is locally approximated by the
non-conforming virtual element space of order $k$ that contains the
subspace of polynomials of degree at most $k$
and is globally non-conforming in the sense specified in
Section~\ref{sect:virtual:element:framework}, see
also~\cite{Ayuso-Lipnikov-Manzini:2016,Cangiani-Manzini-Sutton:2016}.
The pressure is locally approximated by polynomials of degree at most
$k-1$, and is globally discontinuous;
$(ii)$ gradient and divergence are approximated by their projection
onto polynomials of degree $k-1$.
Both projections are computable exactly using only the degrees of
freedom of the VEM.
Therefore, the divergence-free nature of the Stokes velocity is
reproduced in the virtual framework by enforcing the divergence-free
condition on the velocity approximation in a weak sense on each
element.
Moreover, the degree of the polynomials determines the accuracy
(convergence rate) of the VEM;
$(iii)$ the well-posedness of the VEM is ensured through an additional
stabilization term in the discrete weak formulation, which is
computable using only the degrees of freedom of the VEM and is zero
when applied to polynomials;
$(iv)$ the VEM allows for the use of mesh partitionings of $\Omega$
with polygonal elements in 2D or polyhedral elements in 3D of
arbitrary shape provided that a few typical shape regularity
conditions are satisfied. The formulation of the method is the same in
2D and 3D and for any cell shape.
It is worth mentioning that this feature follows from the
nonconforming nature of the formulation, since the virtual conforming
formulation, which also holds for general meshes, has to be
constructed hierarchically in the space dimensions.

A number of relevant numerical approaches for the Stokes problem have been 
proposed in recent years.
The VEM framework has already been applied to the streamline formulation
of the Stokes equation in~\cite{Antonietti-BeiraodaVeiga-Mora-Verani:2014},
a pseudo-stress velocity formulation can be found in~\cite{Caceres-Gatica:2016},
and a divergence free virtual approach can be found
in~\cite{BeiraodaVeiga-Lovadina-Vacca:2016}. 
Among the recent developments in discontinuous Galerkin methods, it is
worth mentioning~\cite{Burman-Stamm:2010}, the Hybridized
Discontinuous
Galerkin~\cite{Cockburn-Shi:2014,Cockburn-Nguyen-Peraire:2010}, and,
concerning polygonal meshes, the Hybrid High Order
method~\cite{Aghili-Boyaval-DiPietro:2014} and the Weak Galerkin
method~\cite{Wang-Ye:2016}.
A comparison between these different approaches is surely worth of
investigation and will be the subject of further investigations.

The paper is organized as follows.
In Section~\ref{sect:stokes:problem} we state the steady Stokes
problem in weak form and introduce the VEM as a Galerkin method.
In Section~\ref{sect:virtual:element:framework} we review the
non-conforming virtual element framework used to approximate the
velocity field, while implementation details can be found in 
Section~\ref{sec:implementation}.
In Section~\ref{sect:error:analysis} we prove the well-posedness and
convergence of the VEM and we derive the error estimates for the
velocity and pressure approximation.
In Section~\ref{sect:numerical:experiments} we numerically assess the
performance of the VEM by solving a set of representative problems.
In Section~\ref{sect:conclusion} we offer our final remarks and
conclusions.

\section{Continuous Stokes problem and discrete formulation}
\label{sect:stokes:problem}
The primary velocity-pressure formulation of the Stokes
problem~\eqref{eq:Stokes:PDE:a}-\eqref{eq:Stokes:PDE:c} takes the
variational form:

\medskip
\emph{Find $\Bu\in\big[H^1(\Omega)\big]^d$ with $\Bu=\Bg\in\big[H^{\frac{1}{2}}(\Gamma)\big]^d$ on $\Gamma$
  and $p\in L^2(\Omega)/\R$ such that for $\Bf\in\big[L^2(\Omega)\big]^d$ it holds:}
\begin{align}
  \as(\Bu,\Bv) + \bs(\Bv, p) &= (\Bf,\Bv) \phantom{0}\qquad\forall\Bv\in\big[H^1_0(\Omega)\big]^d, \label{eq:Stokes:PDE:weak:a}\\[0.5em]
  \bs(\Bu,q)                 &= 0 \phantom{(\Bf,\Bv)}\qquad\forall q \in L^2(\Omega)/\R,           \label{eq:Stokes:PDE:weak:b}
\end{align}
where $\Gamma$ is the boundary of $\Omega$ and the bilinear forms $a$
and $b$ are defined by:
\begin{align}
  \as(\Bu,\Bv) = \int_\Omega \D\Bv:\D\Bu\,\dx,
  \qquad
  \bs(\Bv,q  ) = -\int_\Omega q\,\Div\Bv\,\dx.
  \label{eq:global:bilinear:form}
\end{align}
The well-posedness
of~\eqref{eq:Stokes:PDE:weak:a}-\eqref{eq:Stokes:PDE:weak:b} follows
from the coercivity of the form $\as$ on the kernel of the form $\bs$
and the inf-sup condition~\cite{Boffi-Brezzi-Fortin:2013}.

In~\eqref{eq:Stokes:PDE:weak:a}-\eqref{eq:Stokes:PDE:weak:b} and
throughout the paper we use the standard definitions and notation of
Sobolev spaces, inner products, seminorms and norms.
In particular, if $D$ is an open bounded domain with Lipschitz
boundary in $\R^d$ for $d=2,3$ and $m$ a non-negative integer,
$H^m(D)$ denotes the standard Sobolev space of order $m$;
$(\cdot,\cdot)_{m,D}$ is the associated inner product;
$\norm{\,\cdot\,}_{m,D}$ and $\abs{\,\cdot\,}_{m,D}$ are the induced
norm and seminorm, respectively.
When $D=\Omega$ as
in~\eqref{eq:Stokes:PDE:weak:a}-\eqref{eq:Stokes:PDE:weak:b} we drop
the subscripted symbol.

Let $k\geq 1$ be a fixed integer.
A Virtual Element Method of order $k$ will be defined by two finite
dimensional functional spaces $\VVhk$ and $\Phi^{k-1}_h$ of discrete
trial velocity and pressure fields and bilinear forms
$\ah:\VVhk\times\VVhk\rightarrow\R$ and
$\bh:\VVhk\times\Phk\rightarrow\R$ discrete counterparts of $\as$ and
$\bs$, respectively.
Precise definition of the functional spaces $\VVhk$ and $\Phi^k_h$ and
the construction of the bilinear forms $\ah$ and $\bh$ will be the
focus of most of the remainder of this paper.
For the moment, we only anticipate that we shall not assume the
inclusion $\VVhk\subset\big[H^1(\Omega)\big]^d$ as our main goal is
the development of a non-conforming approximation.
Moreover, let $\Bgh$ be a suitable piecewise polynomial approximation
of $\Bg$ on the mesh partitioning of $\Gamma$.
The precise definition is given at the end of
section~\ref{sect:non-conforming:virtual:element:space}.
The virtual element formulation for the approximate solution
of~\eqref{eq:Stokes:PDE:weak:a}-\eqref{eq:Stokes:PDE:weak:b} reads as:

\medskip
\emph{Find $(\Buh,\ph)\in\VVhkg\times\Phk$ such that}
\begin{align}
  \ah(\Buh,\Bvh) + \bh(\Bvh, \ph) &= (\Bfh,\Bvh) \phantom{0}\qquad\forall\Bvh\in\VVhkz, \label{eq:Stokes:dPDE:weak:a}\\[0.5em]
  \bh(\Buh,\qh)                   &= 0 \phantom{(\Bfh,\Bvh)}\qquad\forall \qh \in\Phk, \label{eq:Stokes:dPDE:weak:b}
\end{align}
with $\VVhkg=\{\Bvh\in\VVhk:\Bvh|_\Gamma=\Bgh\}$ and
$\VVhkz=\{\Bvh\in\VVhk:\Bvh|_\Gamma=0\}$, respectively.
The vector field $\Bfh$ in the right-hand side integral
of~\eqref{eq:Stokes:dPDE:weak:a} is a suitable approximation of the
vector field $\Bf$.
The well-posedness of
problem~\eqref{eq:Stokes:dPDE:weak:a}-\eqref{eq:Stokes:dPDE:weak:b}
will follow from a discrete inf-sup condition which shall be
established under suitable coercivity and stability properties
introduced in the following section.

\section{Virtual element framework}
\label{sect:virtual:element:framework}

\subsection{Mesh regularity and polynomial approximation}
\label{subsec:mesh:jumps:polynomials}
To ease the exposition, we assume that $\Omega$ is a polygonal domain
for $d=2$ and a polyhedral domain for $d=3$.
For any fixed $h>0$ we have a finite decomposition (the mesh) $\Th$ of
the domain $\Omega$ into non-overlapping \emph{simple
  polygonal/polyhedral elements} with maximum size $h$.
The adjective ``simple'' refers to the fact that the boundary of each
element in the decomposition must be non-intersecting.
Moreover, the boundary $\partial E$ of element $E$ is made of a
uniformly bounded number of interfaces (edges/faces), which are either
part of the boundary of $\Omega$, or shared with another element of
the decomposition.
The definition of \emph{simple polygons} and \emph{simple polyhedra}
is general enough to include, for instance, elements with consecutive
co-planar edges/faces, such as those typical of locally refined meshes
with hanging nodes and non-convex elements.

Below, we use $s$ to denote a $d-1$ dimensional mesh interface (either
an edge when $d=2$ or a face when $d=3$), $\abs{s}$ to denote its
length, $\Bn_{s}$ to denote its unit normal vector with orientation
fixed once and for all, and $\Edges$ to denote the set of all such
mesh interfaces in $\Th$.
When referring to the boundary of a specific element $E$ (with $\nu_E$
edges/faces) we use the notation $s\in\d E$ and $\Bn_s$ will have the
outward orientation.

\medskip
\noindent
\begin{assum}[Mesh regularity]
  \label{assum:mesh:regularity}
  We assume that there exists a constant $\rho>0$ such that:
  \begin{itemize}
  \item for every element $E$ of $\Th$ and every interface $s\in E$,
    it holds that $h_s\geq\rho h_E$;
  \item every element $E$ of $\Th$ is star-shaped with respect to a
    ball of radius $\rho h_E$;
  \item for $d=3$, every face $s$ of the mesh is star-shaped with
    respect to a ball of radius $\rho h_s$.
  \end{itemize}
\end{assum}

\medskip
\noindent
If $s$ is an internal edge/face of $\Th$, then, there exist two
elements $E^+$ and $E^-$ such that $s\subset\d E^+\cap\d E^-$.
Consider a scalar function $v$ defined on $\Omega$.
We denote by $v^{\pm}$ the trace of $v_{|_{E^{\pm}}}$ on $s$ from
within $E^{\pm}$ and by $\Bn_s^{\pm}$ the unit vector orthogonal to
$s$ and pointing out of $E^{\pm}$.
Then, the \emph{jump of the scalar function $v$ across $s$} is defined
as $\jump{v}:=v^+\Bn_{s}^{+}+v^-\Bn_{s}^{-}$.
If, on the other hand, $s$ is on the domain boundary $\Gamma$, then
$\jump{v}:=v\Bn_{s}$, with $v$ representing the trace of $v$ from
within the element $E$ having $s$ as an interface and $\Bn_{s}$ is the
unit vector orthogonal to $s$ and pointing out of $\Omega$.
Similarly, the jump of the vector quantity $\Bvh$ at the internal
interface $s$ is given by
\begin{align}
  \jump{\Bvh} = (\Bn_s^+\cdot\Bvh^++\Bn_s^-\cdot\Bvh^-),
  \label{eq:jump:vector:def}
\end{align}
and for the tensor quantity $\Bn\Bvh$ we may consider the jump
operator $\jumpt{\,\cdot\,}$ that is such that
\begin{align}
  \jumpt{\Bvh}:{\bm\tau}=
  \big(\Bn_s^+\cdot{\bm\tau}\big)\cdot\Bvh^++\big(\Bn_s^-\cdot{\bm\tau}\big)\cdot\Bvh^-
  \label{eq:jump:tensor:def}
\end{align}
for every properly sized tensor quantity ${\bm\tau}$.

\medskip
\noindent
We denote by $\PiE{l}\,:\,L^2(E)\to\PE{l}$ for $l\geq 0$ the
$L^2(E)$-orthogonal projection onto the polynomial space $\PE{l}$,
defined for any function $v\in L^2(E)$ as the unique solution of the
problem:
\begin{align}
  \label{eq:l2ProjDefn}
  (\PiE{l}(v),q)_{E} = (v,q)_{E} \quad \forall q\in\PE{l}.
\end{align}
For vector fields, i.e., $\Bv\in\big[L^2(E)\big]^d$,
definition~\eqref{eq:l2ProjDefn} is applied component-wise, thus
giving the vector of polynomials $\PiE{l}(\Bv)\in\big[\PE{l}\big]^d$.
The well known approximation property of $L^2(E)$-orthogonal
projection are summarised in the following theorem.

\medskip
\begin{thm}[Approximation using polynomials]
  \label{thm:polynomialApproximation}
  Under Assumption~\ref{assum:mesh:regularity}, the two following
  propositions hold true.
  \begin{enumerate}
  \item[$(i)$] Let $E\in\Th$ and let
    $\PiE{l}\,:\,L^2(E)\rightarrow\PE{l}$, for $l\ge 0$, denote the
    $L^2(E)$-orthogonal projection onto the polynomial space $\PE{l}$.
    Then, for any $w\in H^{\reg}(E)$, with $1\leq \reg\leq l+1$, it
    holds that
    \begin{equation*}
      \norm{w - \PiE{l}(w)}_{0,E} + h_{E} \abs{w - \PiE{l}(w)}_{1,E} \leq Ch_{E}^\reg \abs{w}_{\reg,E}.
    \end{equation*}
  \item[$(ii)$] Let $s$ be an interface shared by $E^{+},E^{-}\in\Th$
    and let $\PiE{l}\,:\,L^2(s)\rightarrow\mathcal{P}_l(s)$, for $l\ge
    0$, denote the $L^2(s)$-orthogonal projector onto the polynomial
    space $\mathcal{P}_l(s)$.
    Then, for every $w\in H^{\reg}(E^{+}\cup E^{-})$, with
    $1\leq\reg\leq l+1$, it holds
    \begin{equation*}
      \big|w - \Pis{l}(w)\big|_{0,s}
      +
      h_s \big|w - \Pis{l}(w)\big|_{1,s}
      \leq
      Ch_{s}^{\reg-1/2} \norm{w}_{\reg,E^{+}\cup E^{-}}.
    \end{equation*}
  \end{enumerate}
  In both instances $(i)$ and $(ii)$, the positive constant $C$
  depends only on the polynomial degree $l$ and the mesh regularity.
\end{thm}
\BEGINPROOF
This theorem can be proven using the theory
in~\cite{Brenner-Scott:1991} for star-shaped domains and its extension
to more general shaped elements presented in,
e.g.,~\cite{Dupont-Scott:1980}.
\ENDPROOF

\subsection{Discrete pressure space}
As discrete trial space for pressures we use the standard space of
piecewise polynomials of degree up to $k-1$ with respect to the domain
partition $\Th$:
\begin{align*}
  \Phk
  := \left\{\,\qh\in L^2(\Omega)\slash{\R}\,|\,{\qh}_{|E}\in\PE{k-1}\,\forall E\in\Th\,\right\}.
\end{align*}
The local degrees of freedom of $\Phk$ for a pentagonal cell and the
polynomial orders $k=1,\ldots,4$ are illustrated by the right
sub-panels in Fig.~\ref{fig:DoF:2D}.
Note that by definition all functions in $\Phk$ are with global zero
mean.
Therefore, if $N_{d,k-1}^{E}$ is the dimension of $\PE{k-1}$, the total
number of degrees of freedom that are required for the pressure
approximation is equal to $\sum_{E\in\Th}N_{d,k-1}^{E}-1$.

\subsection{Scalar non-conforming virtual element space}
\label{sect:scalar:non-conforming:virtual:element:space}
The \emph{scalar non-conforming virtual element space of order $k\geq
  1$ on the element $E$} is defined for $d=2,3$
as~\cite{Ayuso-Lipnikov-Manzini:2016,Cangiani-Manzini-Sutton:2016}
\begin{align}
  \Vhk(E) = \Big\{\,
  v\in H^1(E)\,|\,\Delta v\in \PE{k-2},\,
  \Bn_{s}\cdot\D v\in \PS{k-1}(s)\,\,
  \forall s\in\d E\,\Big\},
  \label{eq:def:WEs}
\end{align}
with the usual convention that $\PE{-1}=\{0\}$.
The virtual element space $\VhE$ contains the space $\PE{k}$ of
polynomials of degree up to $k$ on $E$.
The complement $\Vhk(E)\backslash{\PE{k}}$ is made up of functions
that are deemed expensive to evaluate, although they can be
represented in a discrete form through their degrees of freedom.
The choice of the degrees of freedom is crucial to ensure that it is
possible to define the bilinear forms
in~\eqref{eq:Stokes:dPDE:weak:a}-\eqref{eq:Stokes:dPDE:weak:b} that
are computable just using the degrees of freedom and the polynomial
component of space $\Vhk(E)$.
In practice, the discrete representation is sufficient for the construction of the method.
To characterize the degrees of freedom of the functions in $\Vhk(E)$,
we first introduce an appropriately scaled basis for $\PE{k}$.
Denote by $\cM^\star_{l}(E)$, $l\in\mathbb{N}$, the set of
\emph{scaled monomials}
\begin{equation*}
  \cM^\star_{l}(E) :=
  \left\{\,\Big(\frac{\Bx - \Bx_{E}}{h_{E}}\Big)^{\alpha},\,\abs{\alpha}=l\,\right\} ,
\end{equation*}
where $\alpha$ is a multi-index and $\Bx_{E}$ the center of gravity of
$E$.
Furthermore, we define
$\cM_{k}(E):=\bigcup_{l\leq k}\cM^\star_{l}(E)=:\{m_\alpha\}_{\alpha=1}^{N_{d,k}}$,
a basis of the polynomial space $\PE{k}$ whose size is $N_{d,k}$.
Bases for polynomial spaces defined on an interface $s$ can be
similarly constructed; the same notation will be used.

The degrees of freedom for the scalar non-conforming space $\Vhk(E)$
are~\cite{Ayuso-Lipnikov-Manzini:2016,Cangiani-Manzini-Sutton:2016}:
\begin{itemize}
\item for $k\geq 1$, the moments of degree $(k-1)$ on each edge/face
  $s\in\d E$:
  \begin{equation}\label{eq:edge:moments}
    \mu_s^{\alpha}(\vh)
    :=
    \frac{1}{|s|}
    \int_s \vh\,\ma\,\ds,\qquad\ma\in\cM_{k-1}(s);
  \end{equation}
\item for $k\geq 2$, the moments of degree $(k-2)$ inside the element
  $E$:
  \begin{equation}\label{eq:cell:moments}
    \mu_E^{\alpha}(\vh)
    :=
    \frac{1}{|E|}
    \int_E\vh\,\ma\,\dx,\qquad\ma\in\cM_{k-2}(E).
  \end{equation}
\end{itemize}
The degrees of freedom for a pentagonal cell and the polynomial orders
$k=1,\ldots,4$ are illustrated by the left sub-panels in
Fig.~\ref{fig:DoF:2D}.

A counting argument shows that the cardinality of the above sets of
degrees of freedom is $N_E = \nu_E N_{d-1,k-1} + N_{d,k-2}$, where we
recall that $\nu_E$ denotes the number of edges/faces of element $E$.
Moreover, they are unisolvent
in~$\Vhk(E)$~\cite{Ayuso-Lipnikov-Manzini:2016}.
Indeed, if all the degrees of freedom of $\vh$ are zero we find that
\begin{equation}\label{eq:unisolvency:proof}
  \norm{\D\vh}_{0,E}^2
  = (\D\vh,\D\vh)_E
  = -(\vh,\Delta\vh)_E + \sum_{s\in\d E} (\vh,\Bn_s\cdot\D\vh)_s
  = 0.
\end{equation}
To prove~\eqref{eq:unisolvency:proof}, note that for $k=1$ it holds
that $\Delta\vh=0$; for $k>1$ we have that $\Delta\vh\in\PE{k-2}$ and
the first term on the right of~\eqref{eq:unisolvency:proof} is a
linear combination of the internal degrees of freedom of $\vh$, which
are zero by hypothesis.
The second term on the right is also zero because it is a linear
combination of the edge/face degrees of freedom of $\vh$, which are
zero by hypothesis.
From $\norm{\D\vh}_{0,E}=0$ it follows that $\vh$ is constant on $E$
and it must be zero since its value is equal to its zero-th order
moment, which is zero by hypothesis.

\medskip The definition of the Virtual Element Method relies on the
availability of elemental projection operators.
The non-conforming VEM for elliptic problems introduced
in~\cite{Ayuso-Lipnikov-Manzini:2016} is based on the Ritz-Galerkin
projection operator $\Pi^{\nabla}_{k}:H^1(E)\rightarrow \PE{k}$ that
for $v\in H^1(E)$ gives $\Pi^{\nabla}_{k}v$ as the solution of the
problem
\begin{equation*}
  ( \nabla(v - \Pi^{\nabla}_{k} v), \nabla\ma)_{E} = 0 \qquad\forall\ma\in\cM_{k}(E),
\end{equation*}
together with the condition
\begin{equation*}
  \displaystyle \int_{\partial E} (v - \Pi^{\nabla}_{k} v) \ds = 0 \quad \text{if $k=1$,}\qquad
  \int_{E} (v - \Pi^{\nabla}_{k} v) \dx = 0 \quad \text{if $k\ge 2$}.
\end{equation*}
This is shown in~\cite{Ayuso-Lipnikov-Manzini:2016} to be computable
for any $\vh\in\VhE$ using only the degrees of
freedom~\eqref{eq:edge:moments} and~\eqref{eq:cell:moments}.
Furthermore, we shall prove in the following section that the
$L^2$-projector $\PiE{k-1}$ of
Theorem~\ref{thm:polynomialApproximation} is also computable when
applied to first order derivatives of virtual functions.
Together, these projectors will permit us to define a virtual
formulation for the Stokes problem.

\begin{remark}
  Instead, we note that an $L^2$-projection $\PiE{k}(\vh)$ onto
  $\PE{k}$ is not available.
  In view of definition~\eqref{eq:l2ProjDefn}, to compute the
  $L^2$-projection we would need the solution of the finite
  dimensional variational problem: \emph{find $\PiE{k}\vh\in\PE{k}$
    such that}
  \begin{equation}
    \label{eq:computingPovh}
    (\PiE{k}(\vh),\ma)_{E} = (\vh,\ma)_{E} \qquad \forall\ma\in\cM_{k}(E),
  \end{equation}
  which needs the internal moments of $\vh$ up to order $k$.
  As the corresponding degrees of freedom are available only up to
  $k-2$, we need to resort to the strategy originally devised
  in~\cite{Ahmad-Alsaedi-Brezzi-Marini-Russo:2013} for the conforming
  virtual element spaces and extended to the non-conforming space
  by~\cite{Cangiani-Manzini-Sutton:2016}.
  The availability of the $L^2$-projection becomes essential when
  low-order terms are present.
  Since in the present work we do not need this projection, we will
  not consider this issue anymore.
  \ENDPROOF
\end{remark}

\medskip
The global scalar non-conforming virtual element space is defined as a
finite dimensional subspace of the non-conforming Sobolev space
$H^{1,\text{nc}}_k(\Th)$.
The latter is a subspace of the broken Sobolev space
\begin{equation*}
  H^1(\Th):=\big\{v\in L^2(\Omega)\,|\,v_{|E}\in H^1(E),\,\,\forall E\in\Th\big\}
\end{equation*}
and, for $k\geq 1$, is given by
\begin{equation*}
  H^{1,\text{nc}}_{k}(\Th) = \left\{
    v\in H^1(\Th)\,|\,\int_s\jump{v}\cdot\Bn_{s}\,q\,\ds=0
    \quad\forall q\in\PS{k-1}(s),\,\forall s\in\Edges
  \right\},
\end{equation*}
where $\jump{v}$ is the \emph{jump of $v$} across the mesh interface
$s\in\Edges$ defined in
subsection~\ref{subsec:mesh:jumps:polynomials}.
The \emph{global scalar non-conforming virtual element space of order
  $k\geq 1$} is now given by
\begin{equation*}
  \Vhk :=
  \left\{\,\vh\in H^{1,\text{nc}}_{k}(\Th)\,|\,{\vh}_{|E}\in\Vhk(E)\,\,\forall E\in\Th\,\right\}.
\end{equation*}
The degrees of freedom of $\Vhk$ are the edge/face
moments~\eqref{eq:edge:moments} and the internal
moments~\eqref{eq:cell:moments}, and the size of $\Vhk$ is clearly
given by $N_{\Th}\,N_{d,k-2} + N_{d-1,k-1}\,N_{\Edges}$, where
$N_{\Th}$ is the number of cells in $\Th$ and $N_{\Edges}$ the number
of edges/faces in $\Edges$. 
Note that the edge/face moments are the same for the two mesh cells
sharing a given internal edge/face.
Therefore, the weak continuity condition on the jumps in the
definition of $H^{1,\text{nc}}_{k}(\Th)$ is automatically satisfied.

\subsection{Discrete velocity space}
\label{sect:non-conforming:virtual:element:space}

\begin{figure}[t!]
  \[
  \qquad
  \includegraphics[width=.90\textwidth]{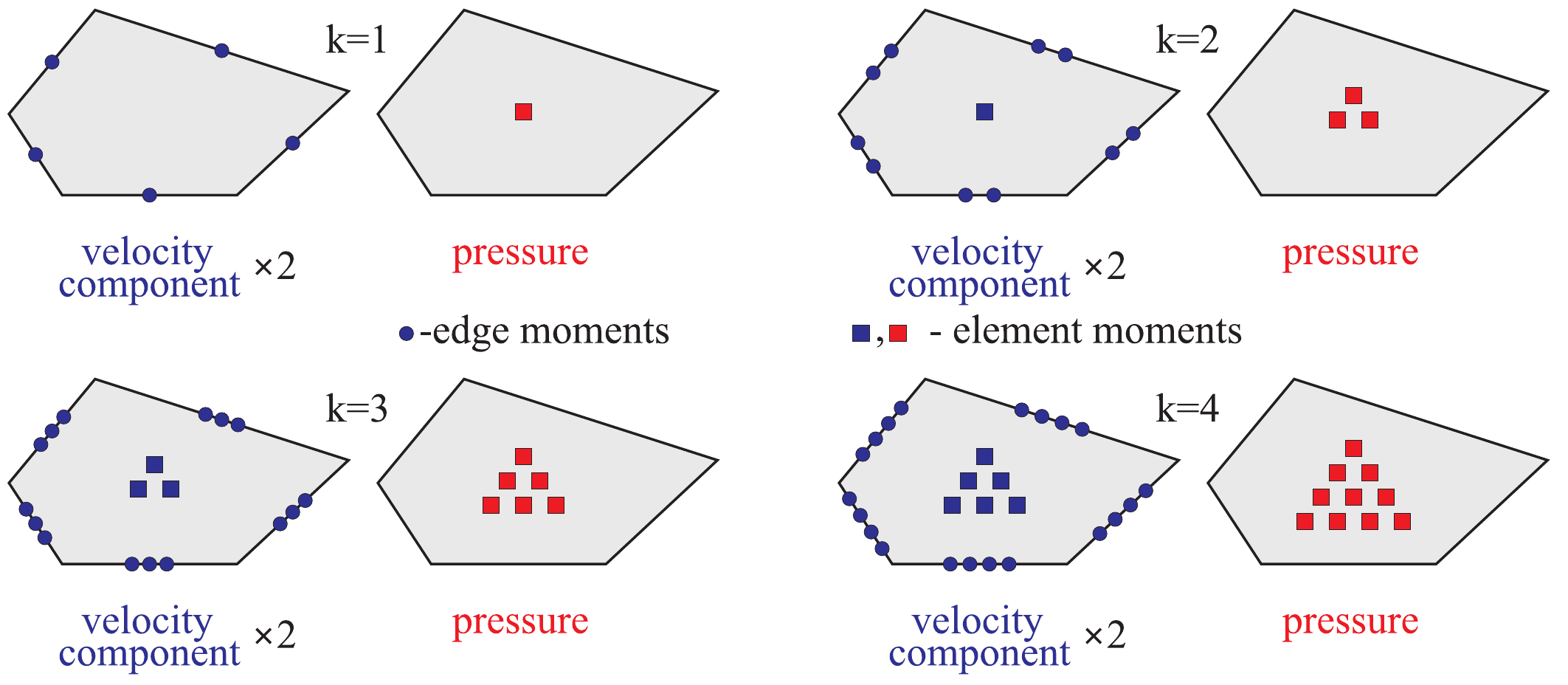}
  \]
  \vspace{-1cm}
  \caption{Illustration of the degrees of freedom for the velocity and
    pressure solving the two-dimensional Stokes
    problem.}
  \label{fig:DoF:2D}
  \medskip
\end{figure}

The discrete velocity space is defined by using the scalar
non-conforming virtual element space of the previous subsection for
each component.
Hence,
\begin{equation*}
\VVhk(E)=[\Vhk(E)\big]^d\qquad\forall E\in\Th,
\end{equation*}
and, similarly, $\VVhk=[\Vhk\big]^d$.
The degrees of freedom of $\VVhk$ are those inherited from each
component, and are illustrated by the left sub-panels in
Fig.~\ref{fig:DoF:2D} for a pentagonal cell and the polynomial orders
$k=1,\ldots,4$.
This choice of degrees of freedom ensures that the discrete
differential operators
\begin{equation}
  \label{eq:projectionsToCompute}
  \PiE{k-1}\circ\nabla
  \quad\textrm{and}\quad
  \PiE{k-1}\circ\Div,
\end{equation}
when applied to a vector field $\Bvh$ in $\VVhk(E)$ are computable
using only the degrees of freedom of $\Bvh$ on $E$.
The projection $\PiE{k-1}\D\Bvh$ is defined as
\begin{align*}
  (\PiE{k-1}(\D\Bvh),\Bma)_{E} = (\D\Bvh,\Bma)_{E}\qquad\forall\Bma\in\big[\cM_{k-1}(E)\big]^{d}.
\end{align*}
Integration by parts yields:
\begin{align*}
  (\D\Bvh,\Bma)_{E} = -(\Bvh,\D\Bma)_{E} + \sum_{s\in\d E}(\Bn_s\cdot\Bvh,\Bma)_{s}.
\end{align*}
Likewise, the projection $\PiE{k-1}(\Div\Bvh)$ is defined as
\begin{align*}
  (\PiE{k-1}(\Div\Bvh),\ma)_{E} = (\Div\Bvh,\ma)_{E}\qquad\forall\ma\in\cM_{k-1}(E).
\end{align*}
Integration by parts yields:
\begin{align*}
  (\Div\Bvh,\ma)_{E} = -(\Bvh,\D\ma)_{E} + \sum_{s\in\d E}(\Bn_s\cdot\Bvh,\ma)_{s}.
\end{align*}
The components of $\Bma\in\big[\cM_{k-1}(E)\big]^d$ and
$\ma\in\cM_{k-1}(E)$ are polynomials of degree at most $k-2$ in $E$
and, when restricted to each $s\in\d E$, are polynomials of degree at
most $k-1$.
Therefore, the right-hand side of the last equation above is
computable using only the internal and the edge/face degrees of
freedom of the (scalar) components $\Bvh$.

\medskip
Finally, in the virtual element formulation
\eqref{eq:Stokes:dPDE:weak:a}-\eqref{eq:Stokes:dPDE:weak:b}, we use
the virtual element space $\VVhkg$, whose definition requires the
boundary function $\Bgh$.
This function is such that ${\Bgh}_{|s}$ on every edge/face
$s\in\Gamma$ is the $L^2$-orthogonal projection of $\Bg$ on the
polynomial space $\PS{k-1}(s)$.
  
\medskip
\subsection{Approximation of the bilinear forms $\as(\cdot,\cdot)$ and $\bs(\cdot,\cdot)$}
\label{subsec:VEM:bilinear:forms}
In view of the following analysis, it is useful to extend the
definition of the continuous bilinear forms $\as$ and $\bs$, to the
whole of $H^1(\Th)$ as a sum of the elemental contributions $\aE$ and
$\bE$,
\begin{align*}
  \as(\Bu,\Bv)=\sum_{E\in\Th}\aE(\Bu,\Bv)
  \quad\textrm{and}\quad
  \bs(\Bu,q)=\sum_{E\in\Th}\bE(\Bu,q)
  \qquad \forall\Bu,\Bv\in H^1(\Th),\,q\in L^2(\R),
\end{align*}
where $\aE$ and $\bE$ are defined by restricting the integrals
in~\eqref{eq:global:bilinear:form} to $E$.

We define the approximate bilinear forms $\ah$ and $\bh$ used
in~\eqref{eq:Stokes:dPDE:weak:a}-\eqref{eq:Stokes:dPDE:weak:b} by
splitting them into local contributions
\begin{align*}
  \ah(\Buh,\Bvh):=\sum_{E\in\Th}\ahE(\Buh,\Bvh)
  \quad\textrm{and}\quad
  \bh(\Buh,\qh):=\sum_{E\in\Th}\bhE(\Buh,\qh),
\end{align*}
for any $\Buh,\Bvh\in\VVhk$ and $\qh\in\Phk$, where $\ahE$ and $\bhE$
are bilinear forms on $\VVhk(E)\times\VVhk(E)$ and
$\VVhk(E)\times\PE{k-1}$, respectively.

The former bilinear form is defined by:
\begin{align}
  \label{eq:ahE}
  \ahE(\Buh,\Bvh) := \int_E\PiE{k-1}(\D\Buh):\PiE{k-1}(\D\Bvh)\,\dx +
  \shE\left(\,(I-\Pi^{\nabla}_{k})\Buh, (I-\Pi^{\nabla}_{k})\Bvh\,\right),
\end{align}
where $\Pi^{\nabla}_{k}$ represents the Ritz-Galerkin projection
operator introduced in
Section~\ref{sect:scalar:non-conforming:virtual:element:space} applied
component-wisely.
The term $\shE$ is the \emph{VEM stabilization} term,
cf.~\cite{BeiraodaVeiga-Brezzi-Cangiani-Manzini-Marini-Russo:2013}.
This  can be {\em any} symmetric and coercive bilinear form
satisfying
\begin{equation}\label{eq:local:stabterm}
  c_*\as(\Bvh,\Bvh)\leq\shE(\Bvh,\Bvh)\leq c^*\as(\Bvh,\Bvh)
  \qquad\text{for all}\quad
  \Bvh\in\VVhk(E)\backslash{\big[\PS{k}(E)]^d},
\end{equation}
for two positive constants $c_*$ and $c^*$ independent of $h$ and the
mesh element $E$. 

Following~\cite{BeiraodaVeiga-Brezzi-Cangiani-Manzini-Marini-Russo:2013,Ayuso-Lipnikov-Manzini:2016},
in all computations presented in
Section~\ref{sect:numerical:experiments} we used the choice
\begin{equation}\label{eq:stabterm}
  \shE\left(\,(I-\Pi^{\nabla}_{k})\Buh, (I-\Pi^{\nabla}_{k})\Bvh\,\right) 
  = \sum_{i=1}^{N_E}\bm{\chi}_i\left((I-\Pi^{\nabla}_{k})\Buh\right)\cdot  \bm{\chi}_i\left((I-\Pi^{\nabla}_{k})\Bvh\right), 
\end{equation}
where $\bm{\chi}_i$ is the vector-valued linear operator that
associates any virtual function $\bm{\phi}$ with the vector of its
$i$-th local degrees of freedom $\bm{\chi}_i(\bm{\phi})\in \R^d$.

\begin{remark}
  Again following
  \cite{BeiraodaVeiga-Brezzi-Cangiani-Manzini-Marini-Russo:2013,Ayuso-Lipnikov-Manzini:2016},
  we may have defined the consistency term on the right-hand side
  of~\eqref{eq:ahE} as $\int_E \D\Pi^{\nabla}_{k}\Buh :
  \D\Pi^{\nabla}_{k}\Bvh \,\dx$, but the approach used
  in~\eqref{eq:ahE} is more suitable to generalisations to problems
  with non-constant coefficients, see,
  e.g.,~\cite{Cangiani-Manzini-Sutton:2016}.
  Furthermore, the availability of the projection $\PiE{k-1}(\D\Bvh)$
  for all $\Bvh\in\VVhk$ proven in the previous section is of its own
  interest.
  \ENDPROOF
\end{remark}

The second bilinear form is defined by:
\begin{align}
  \bhE(\Bvh,\qh) := \int_E\qh\PiE{k-1}(\Div\Bvh)\,\dx.
  \label{eq:bhE:def}
\end{align}

\medskip
\begin{remark}
  \label{rem:bhE:bE}
  We note that
  $\bhE(\cdot,\qh)=\bE(\cdot,\qh)$ in $H^{1,\text{nc}}_{k}(\Th)$ for any $\qh\in\Phk$.
  \ENDPROOF
\end{remark}

\medskip
Definition~\eqref{eq:ahE} with the VEM stabilisation term
satisfying~\eqref{eq:local:stabterm} guarantees that the following
\emph{polynomial consistency} and \emph{stability} properties are
satisfied by the bilinear form $\ahE$.

\medskip
\begin{lem}[Consistency and Stability]
  ~\\ \vspace{-0.5\baselineskip}
  \begin{enumerate}
  \item[$(i)$] \emph{Polynomial consistency}: If $\Buh$ or $\Bvh$, or
    both, belong to $\big[\PE{k}\big]^d$, the bilinear form $\ahE$
    satisfies
    \begin{align}
      \ahE(\Buh,\Bvh) = \aE(\Buh,\Bvh).
      \label{eq:pcons:conds}
    \end{align}
  \item[$(ii)$] \emph{Stability}: There exist two positive constants
    $\alpha_*$ and $\alpha^*$ independent of $h$ and the mesh element
    $E$ such that, for all $\Bvh\in\VVhk(E)$, the bilinear form $\ahE$
    satisfies
    \begin{equation}
      \alpha_*\aE(\Bvh,\Bvh)\leq\ahE(\Bvh,\Bvh)\leq\alpha^*\aE(\Bvh,\Bvh). 
      \label{eq:stab:conds}
    \end{equation}
  \end{enumerate}
\end{lem}
\BEGINPROOF
Property $(i)$ is a straightforward consequence of the fact that the
stabilization term is zero on polynomial vectors.
To prove Property $(ii)$ we first show that $\ahE(\Bvh,\Bvh)=0$
implies that $\Bvh$ is a constant vector, i.e., $\aE$ and $\ahE$ have
the same kernel.
Indeed, consider $\Bvh\in\big[\Vhk(E)\big]^d$ such that
$\ahE(\Bvh,\Bvh)=0$.
We find that
\begin{align*}
  (a)\,\,\int_E\abs{\PiE{k-1}(\D\Bvh)}^2\,\dx = 0
  \quad\textrm{and}\quad
  (b)\,\,\shE\big(\,(I-\PiE{k})\Bvh, (I-\PiE{k})\Bvh\,\big)=0.
\end{align*}
The coercivity of $\shE$ and property $(b)$ imply that
${\Bvh}_{|E}=\PiE{k}(\Bvh)$, i.e., ${\Bvh}_{|E}\in\big[\PE{k}\big]^d$,
and, thus,
${\D\Bvh}_{|E}=\PiE{k-1}(\D\Bvh)\in\big[\PE{k-1}\big]^{d\times d}$.
Property $(a)$ implies that $\PiE{k-1}(\D\Bvh)=0$.
Therefore, it holds that ${\D\Bvh}_{|E}=0$; hence, ${\Bvh}_{|E}$ is a
constant vector.
Now, \eqref{eq:stab:conds} follows from~\eqref{eq:local:stabterm} as
in the scalar case, see~\cite{Cangiani-Manzini-Sutton:2016} for
details.
\ENDPROOF

\vspace{-0.5\baselineskip}
\begin{remark}
  \label{remark:consistency:stability}
  As is usual in the virtual element methodology, in order to satisfy
  conditions \eqref{eq:pcons:conds} and \eqref{eq:stab:conds}, the
  bilinear form $\ahE(\Buh,\Bvh)$ is built as the sum of a ``{\em
    consistency}'' and a ``{\em stabilising}'' term, corresponding to
  the first and the second term in the right-hand side
  of~\eqref{eq:ahE}, respectively.
  The consistency term is exactly computable using only the degrees of
  freedom of $\Buh$ and $\Bvh$ and satisfies the consistency
  condition~\eqref{eq:pcons:conds}.
  However, the consistency term alone does not satisfy the stability
  condition~\eqref{eq:stab:conds} on the whole virtual element space
  due to a rank deficiency of the operator, or equivalently, to the
  existence of a spurious kernel.
  Hence, a \emph{stabilization} term must be added to fix this issue
  and actually remove the spurious kernel.
  This latter term is designed to vanish on the polynomial subspace
  not to affect the consistency property.
  As such, choice of the stabilization term is not
  unique~\cite{Cangiani-Manzini-Sutton:2016}.
  For example, for the closely related MFD method, a number of studies
  investigated the optimal choice of the stabilization term with
  respect to a given criterion (reduction of dispersion effects,
  existence of a discrete maximum/minimum principle),
  cf.~\cite{Gyrya-Lipnikov-Manzini-Svyatskiy:2014,Gyrya-Lipnikov:12,Bokil-Gibson-Gyrya-McGregor:2015}.
  By exploiting the strict relation existing between the MFD method
  and the VEM, these alternative constructions of the stabilization
  term could be optionally considered in the present context.
  \ENDPROOF
\end{remark}

\vspace{-0.5\baselineskip}

\subsection{Mesh-dependent energy norms}
Hereafter, we shall use the energy semi-norm on the broken Sobolev
space $\big[H^1(\Th)\big]^d$:
\begin{align*}
\abs{\Bv}_{1,h}^2 := \sum_{E\in\Th}\abs{\Bv}_{1,E}^2
\quad\textrm{with}\quad
\abs{\Bv}_{1,E}^2 = \aE(\Bv,\Bv).
\end{align*}
A standard application of the results in~\cite{Brenner:2003} shows
that a Poincar\'e inequality holds for the functions in
$H_k^{1,\text{nc}}(\Th)$, $k\geq 1$.
Therefore, the semi-norm $\abs{\cdot}_{1,h}$ is a norm in
$\big[H_k^{1,\text{nc}}(\Th)\big]^d$, and for this reason throughout
the paper we prefer to use the notation
$\bSnorm{\Bv}$ and
$\bSnormE{\Bv}$ instead of
$\abs{\Bv}_{1,h}$ and
$\abs{\Bv}_{1,E}$.
To prove the inf-sup stability, we use the mesh-dependent energy
seminorm on $\VVhk$ given by
\begin{align}
  \TNORM{\Bvh}^2 := \sum_{E\in\Th}\TNORME{\Bvh}^2,
  \label{eq:normh:def:a}
\end{align}
with
\begin{align}
  \TNORME{\Bvh}^2 := \ahE(\Bvh,\Bvh) =
  \int_E\abs{\PiE{k-1}(\D\Bvh)}^2\,\dx +
  \shE\left(\,(I-\PiE{k})\Bvh, (I-\PiE{k})\Bvh\,\right).
  \label{eq:normh:def:b}
\end{align}
Also, we will consider the affine subspace of vector-valued functions
$\VVhkg=\big[\Vhkg\big]^d$ where $\Vhkg$ contains the scalar functions
of $\Vhk$ whose trace on the boundary $\Gamma$ is equal to $g$, i.e.,
\begin{align}
  \Vhkg := \Big\{\vh\in\Vhk\,|\,{\vh}_{|s}=g\textrm{~for~}s\subset\Gamma\Big\},
\end{align}
and the linear subspace $\Vhkz$ obtained for $g=0$.
On $\Vhkz$ we have the following equivalence of norms.

\medskip
\begin{lem}
  \label{lem:normh}
  The seminorm $\TNORM{\,\cdot\,}$ defined
  by~\eqref{eq:normh:def:a}-\eqref{eq:normh:def:b} is a norm on
  $\VVhkz$ and
  \begin{align}
    \sqrt{\alpha_*}\bSnorm{\Bvh}
    \leq\TNORM{\Bvh}\leq
    \sqrt{\alpha^*}\bSnorm{\Bvh}
    \qquad\forall\Bvh\in\VVhkz.
  \label{eq:equiv:norms}
  \end{align}
\end{lem}
\BEGINPROOF
As noted in the proof of the stability
condition~\eqref{eq:stab:conds}, $\TNORME{\Bvh}=0$ implies that
${\Bvh}_{|E}$ is a constant vector, and, from the definition of
$H^{1,\text{nc}}_{k}(\Th)$, it follows that $\Bvh=$ constant.
Finally, from $\Bvh\in\VVhkz$ it follows that $\Bvh=0$.
The equivalence between the two norms in~\eqref{eq:equiv:norms} is a
straightforward consequence of~\eqref{eq:stab:conds}.
\ENDPROOF

\subsection{Approximation of the right-hand side $(\Bf,\cdot)$}
\label{subsec:rhs:approximation}
We approximate the right-hand side term $(\Bf,\,\cdot\,)$ by the
linear functional
\begin{align}
  (\Bfh,\Bvh):=\sum_{E\in\Th}(\Bfh,\Bvh)_E
  \quad\textrm{with}\quad
  {\Bfh}_{|E}:=\PiE{\max(k-2,0)}(\Bf).
  \label{eq:Bfh:def}
\end{align}
Since $\PiE{\max(k-2,0)}(\Bf)$ is a polynomial of degree at most
$k-2$, each local linear functional is bounded and, for $k\geq 2$, it
is computable by using the internal degrees of freedom of $\Bvh$.
However, it is not computable for $k=1$.
Indeed, the computation of $(\Bfh,\Bvh)_E$ above requires the
knowledge of the average value of the components of $\Bvh$ on each
element $E$ and such information is only available for $k\geq 2$.
Therefore, the case $k=1$ deserves a special treatment.
We follow Reference~\cite{Ayuso-Lipnikov-Manzini:2016} and approximate
$\PiE{0}(\Bvh)$ by the average of the $0$-th order moments of $\Bvh$
associated with the edge/face of cell $E$.
Namely we consider
\begin{align*}
  \Bvhb|_{E}=\frac{1}{\nu_E}\sum_{s\in\partial E}\frac{1}{|s|}\int_{s}\Bvh ds,
\end{align*}
and note that $\Bvhb|_{E}$ is a first-order approximation to
$\PiE{0}(\Bvh)=\frac{1}{|E|}\displaystyle\int_{E}\Bvh\dx$, i.e., we
have that
\begin{align*}
\norm{\Bvhb|_{E}-\PiE{0}(\Bvh)}_{0,E}\leq Ch|\Bv|_{1,E}.
\end{align*}
Then, we use $\Bvhb$ to compute $(\Bfh,\Bvh)$ through the
approximation:
\begin{align*}
  \big(\PiE{0}(\Bf),\Bvh\big)_{E}
  =\big(\PiE{0}(\Bf),\PiE{0}(\Bvh)\big)_{E}
  \approx\big(\PiE{0}(\Bf),\Bvhb\big)_{E}.
\end{align*}
Therefore, for $k=1$ we take:
\begin{align*}
  (\Bfh,\Bvh)
  :=\sum_{E\in\Th} (\PiE{0}(\Bf),\Bvhb)_{E}.
\end{align*}
We collect the results for the approximation of the right-hand side
functional $(\Bf,\cdot)$ for $k=1$ and $k>1$ in the following lemma.
The proof follows from a straightforward extension of the scalar case,
which is found in~\cite{Ayuso-Lipnikov-Manzini:2016}, to $d$-sized
vector-valued forcing terms $\Bf$ and for this reason is omitted.
\begin{lem}[Approximation of the right-hand side $(\Bf,\cdot)$]
  \label{lem:rhs:error}
  Let $\reg,k\geq 1$ be integer numbers and consider
  $\Bf\in\big[H^{\reg-1}(\Omega)\big]^d$, $r=\min(k,\reg)$, and
  $(\Bfh,\cdot)\in\big[(\Vhk)'\big]^d$ defined as above.
  There exists a constant $C$ independent of $h$ such that
  \begin{align*}
    \sup_{\Bvh\in\VVhk}\frac{\abs{(\Bfh,\Bvh)-(\Bf,\Bvh)}}{\norm{\Bvh}_{1,h}}
    \leq C h^{r} \norm{\Bf}_{r-1}.
  \end{align*}
\end{lem}
\BEGINPROOF
See~\cite{Ayuso-Lipnikov-Manzini:2016}.
\ENDPROOF

\section{Error analysis}
\label{sect:error:analysis}
The well-posedness of the discrete
problem~\eqref{eq:Stokes:dPDE:weak:a}-\eqref{eq:Stokes:dPDE:weak:b} is
discussed in Section~\ref{subsect:existence:and:uniqueness}.
The non-conformity error is estimated in
Section~\ref{subsect:estimate:non-conformiy:error}.
The convergence analysis is carried out in Sections
\ref{subsect:error:estimates:velocity} and
\ref{subsect:error:estimates:pressure}, where error estimates for the
approximation of the velocity and pressure fields are respectively
derived.
For the convergence analysis, we assume that $\Bg=0$ on $\Gamma$
in~\eqref{eq:Stokes:PDE:c} through the whole section.

Define the \emph{virtual interpolant} of the vector field
$\Bv\in\big[H^{1,\text{nc}}_{k}(\Th)\big]^d$ as the unique vector
field $\BvI\in\VVhk$ whose degrees of freedom are the internal and
edge/face moments of $\Bv$.
Formally,
\begin{itemize}
\item for $k\geq 1$, the degrees of freedom of $\BvI$ associated with
  the edge/face $s\in\Edges$ are given by
  \begin{equation}\label{eq:edge:moments:BvI}
    \mu_s^{\alpha}(\BvI) := \int_s\Bv\,\ma\,\ds,\qquad\ma\in\cM_{k-1}(s);
  \end{equation}
\item for $k\geq 2$, the degrees of freedom of $\BvI$ associated with
  the mesh element $E$ are given by
  \begin{equation}\label{eq:cell:moments:BvI}
    \mu_E^{\alpha}(\BvI) := \int_E\Bv\,\ma\,\dx,\qquad\ma\in\cM_{k-2}(E).
  \end{equation}
\end{itemize}
The above relations must be interpreted component-wise.
The unisolvence of the degrees of freedom implies the uniqueness of
$\BvI$.
Moreover, if $\Bv\in\big[H^1_0(\Omega)\big]^d$ all its moments on each
edge/face $s$ on $\Gamma$ are zero and, consequently, $\BvI$ belongs
to $\VVhkz$.
We also have the following result regarding the approximation of
sufficiently smooth functions by the virtual interpolant, which may be
proven as in~\cite{Ayuso-Lipnikov-Manzini:2016}.

\medskip
\noindent
\begin{thm}[Approximation using virtual element functions]
  \label{thm:spaceApproximation}
  Let $\VVhk$ the non-conforming virtual element space of
  Section~\ref{sect:non-conforming:virtual:element:space} for any
  integer $k\geq 1$, $\reg$ a positive integer such that
  $2\leq\reg\leq k+1$, and $D$ a closed subset of $\Omega$.
  Under Assumption~\ref{assum:mesh:regularity} (mesh regularity), for
  any $\Bv\in H^\reg(D)$, there exists an element $\BvI\in\VVhk$ such
  that
  \begin{equation*}
    \norm{\Bv-\BvI}_{0,D} + h\abs{\Bv - \BvI}_{1,D} \leq Ch^\reg \abs{\Bv}_{\reg,D},
  \end{equation*}
  where $C$ is a positive constant that depends only on the polynomial
  degree $k$ and the mesh regularity constant $\rho$.
\end{thm}

\subsection{Existence and uniqueness of the virtual element solution}
\label{subsect:existence:and:uniqueness}

The main result of this section is the existence and uniqueness of the
virtual element solution $(\Buh,\ph)\in\VVhk\times\Phk$, which is
stated in Theorem~\ref{thm:existence:uniquenessfor:eliptic:problem}.
The proof of this theorem is based on the \emph{inf-sup} property that
is proven in the following lemma by adapting a classical argument.

\medskip
\begin{lem}[Inf-sup]
  \label{lem:inf-sup}
  There exists a strictly positive constant $\beta$ independent of $h$
  such that for every $\qh$ in $\Phk$ there exists a vector $\Bvq$ in
  $\VVhkz$ such that
  \begin{align*}
    \frac{\bh(\Bvq,\qh)}{\TNORM{\Bvq}}\geq \beta\norm{\qh}_{0}.
  \end{align*}
\end{lem}
\BEGINPROOF
From~\cite{Boffi-Brezzi-Fortin:2013} we know that there exists a
strictly positive constant $\widetilde{\beta}$ independent of $h$ such
that for every $q$ in $L^2(\Omega)\slash{\R}$ there exists a vector
$\Bvtq$ in $\big[H^1_0(\Omega)\big]^d$ such that
\begin{align*}
  \frac{\bs(\Bvtq,q)}{\norm{\Bvtq}_{1}}\geq\widetilde{\beta}\norm{q}_{0}.
\end{align*}
We can restrict this inequality to $\Phk\subset L^2(\Omega)\slash{\R}$
and for any $\qh\in\Phk$ consider the corresponding vector $\Bvtq$.
We will prove that
\begin{align*}
  (i)\quad\bh(\BvtqI,\qh) = \bs(\Bvtq,\qh) \quad\forall\qh\in\Phk
  \quad\textrm{and}\quad
  (ii)\quad\TNORM{\BvtqI}\leq \gamma\norm{\Bvtq}_{1},
  \end{align*}
  where $\BvtqI$ is the virtual interpolation of $\Bvtq$ defined in
  the previous section and $\gamma=\alpha^*\slash{\sqrt{\alpha_*}}$.
Such properties easily imply that
\begin{align*}
  \frac{\bh(\BvtqI,\qh)}{\TNORM{\BvtqI}}\geq\frac{\bs(\Bvtq,\qh)}{\gamma\norm{\Bvtq}_{1}}
  \qquad\forall\qh\in\Phk,
\end{align*}
from which the assertion of the lemma follows with $\Bvq=\BvtqI$ and
$\beta=\widetilde{\beta}\slash{\gamma}$.

To prove $(i)$, consider the following development that starts from
the definition of $\bhE$:
\begin{equation*}
  \begin{array}{rll}
    \bhE(\BvtqI,\qh)
    &\displaystyle= \int_E \qh\,\PiE{k-1}(\Div\BvtqI)\dx                                      &\quad\mbox{[use the definition of $\PiE{k-1}$]}\\[0.75em]
    &\displaystyle= \int_E \qh\,\Div\BvtqI\dx                                                &\quad\mbox{[integrate by parts]}\\[0.75em]
    &\displaystyle= -\int_E \D\qh\cdot\BvtqI\dx + \sum_{s\in\d E}\int_s \qh\Bn_s\cdot\BvtqI\ds &\quad\mbox{[use the definition of $\BvtqI$]}\\[0.75em]
    &\displaystyle= -\int_E \D\qh\cdot\Bvtq\dx + \sum_{s\in\d E}\int_s \qh\Bn_s\cdot\Bvtq\ds   &\quad\mbox{[integrate by parts back]}\\[0.75em]
    &\displaystyle= \int_E \qh\,\Div\Bvtq\dx                                                 &\quad\mbox{[use the definition of $\bE(\cdot,\cdot)$]}\\[0.75em]
    &\displaystyle= \bE(\Bvtq,\qh)
  \end{array}
\end{equation*}
Property $(i)$ readily follows from the elemental decomposition of
$\bh$ and $\bs$.

To prove $(ii)$, first note that the stability condition of $\ahE$
implies that
\begin{align}
  \alpha_*\bSnormE{\BvtqI}^2\leq
  \TNORME{\BvtqI}^2
  \leq\alpha^*\bSnormE{\BvtqI}^2.
  \label{eq:stab-vtqI}
\end{align}
Now, since $\BvtqI\in\Vhk(E)$ it holds that $\Delta\BvtqI$ is a vector
of polynomials of degree $k-2$ inside $E$ and $(\Bn_s\cdot\D)\BvtqI$
is a vector of polynomials of degree $k-1$ along each edge $s\in\d E$
and we have that:
\begin{equation*}
  \begin{array}{rll}
    \aE(\BvtqI,\BvtqI)
    &\displaystyle=  \int_E\D\BvtqI:\D\BvtqI\,\dx                                                                   & \quad\mbox{[integrate by parts]}               \\[0.5em]
    &\displaystyle= -\int_E \BvtqI\cdot\Delta\BvtqI\,\dx + \sum_{s\in\d E}\int_s\BvtqI\cdot(\Bn_s\cdot\D)\BvtqI\,\ds  & \quad\mbox{[use the definition of $\BvtqI$]}   \\[0.5em]
    &\displaystyle= -\int_E \Bvtq \cdot\Delta\BvtqI\,\dx + \sum_{s\in\d E}\int_s\Bvtq \cdot(\Bn_s\cdot\D)\BvtqI\,\ds  & \quad\mbox{[integrate by parts back]}          \\[0.5em]
    &\displaystyle=  \int_E\D\Bvtq:\D\BvtqI\,\dx                                                                    & \quad\mbox{[use the Cauchy-Schwarz inequality]}\\[0.5em]
    &\leq \big(\aE(\Bvtq,\Bvtq)\big)^{\frac{1}{2}}\big(\aE(\BvtqI,\BvtqI)\big)^{\frac{1}{2}}                     & \quad\mbox{[use the left inequality of~\eqref{eq:stab-vtqI}]}\\[0.5em]
    &\leq \bSnormE{\Bvtq}\,\frac{1}{\sqrt{\alpha_*}}\TNORME{\BvtqI}.
  \end{array}
\end{equation*}
Since $\bSnormE{\BvtqI}^2=\aE(\BvtqI,\BvtqI)$, the last inequality
and~\eqref{eq:stab-vtqI} implies that
$\TNORME{\BvtqI}\leq\alpha^*\slash{\sqrt{\alpha_*}}\bSnormE{\Bvtq}$
and property $(ii)$ follows from the continuity of $\Bvtq$ in
$\big[H^1_0(\Omega)\big]^d$ by summing over all $E\in\Th$ and setting
$\gamma=\alpha^*\slash{\sqrt{\alpha_*}}$.
\ENDPROOF

\medskip
\noindent
\begin{thm}\label{thm:existence:uniquenessfor:eliptic:problem}
  The solution of the discrete
  problem~\eqref{eq:Stokes:dPDE:weak:a}-\eqref{eq:Stokes:dPDE:weak:b}
  exists and is unique.
\end{thm}
\BEGINPROOF
Remark~\ref{rem:bhE:bE} implies that
$\textsf{ker}(\bh)=\textsf{ker}(\bs)$ in
$\big[H^{1,\text{nc}}_{k}(\Th)\big]^d\times\Phk$, where
\begin{align*}
  \textsf{ker}(\bh) =
  \big\{\Bvh\in\big[H^{1,\text{nc}}_{k}(\Th)\big]^d\,|\,
  \bh(\Bvh,\qh)=0\,\,\forall\qh\in\Phk\big\}.
\end{align*}
In view of the stability condition~\eqref{eq:stab:conds}, the
coercivity of $\ah(\cdot,\cdot)$ on $\textsf{ker}(\bh)$ follows from
the coercivity of $\as(\cdot,\cdot)$ on $\textsf{ker}(\bs)$.
Therefore, existence and uniqueness of the solution follow from the
\emph{inf-sup} property proved in Lemma~\ref{lem:inf-sup},
cf.~\cite{Boffi-Brezzi-Fortin:2013}.
\ENDPROOF

\subsection{Estimate of the non-conformity error}
\label{subsect:estimate:non-conformiy:error}

The non-conformity error is controlled as in
Lemma~\ref{lem:non-conformity:error}.
The proof of this lemma requires a bound on the jumps of $\Bvh$ and
$\Bn\Bvh$ (which, we recall, are defined in~\eqref{eq:jump:vector:def}
and~\eqref{eq:jump:tensor:def}).
This bound is provided by Lemma~\ref{lem:estimates:of:jumps}.

\medskip
\begin{lem}\label{lem:estimates:of:jumps}
  Let $k,\reg\geq 1$, $r=\min(k,\reg)$ be integer numbers.
  Consider $\Bu\in\big[H^{\reg+1}(\Omega)\big]^d$ and $p\in
  L^2(\Omega)\slash{\R}\cap H^{\reg}(\Omega)$.
  Then, under Assumption~\ref{assum:mesh:regularity} (mesh
  regularity), for every $\Bvh\in\VVhkz$ it holds
  \begin{align*}
    \abs{\sum_{s\in\Edges}\int_s\nabla\Bu:\jumpt{\Bvh}\ds\,} +
    \abs{\sum_{s\in\Edges}\int_s p\jump{\Bvh}\ds\,}
    \leq Ch^{r}\big(\norm{u}_{r+1}+\norm{p}_{r}\big)\,\norm{\Bvh}_{1}.
  \end{align*}
\end{lem}
\BEGINPROOF
%
%
Since $\Bvh\in\VVhk\subset\big[H^{1,\text{nc}}_{k}(\Th)\big]^d$, the
jump of the components of $\Bvh$ on every interface $s$ is orthogonal
to the polynomial functions of degree $k-1$ defined along $s$.
Thus, it holds that
\begin{align*}
  \sum_{s\in\Edges} \int_s p\jump{\Bvh}\,\ds
  =
  \sum_{s\in\Edges}
  \int_s \big(p -\Pis{k-1}(p)\big)\,\big(\jump{\Bvh}-\Pis{k-1}\jump{\Bvh}\big)\ds,
\end{align*}
for the pressure and
\begin{align*}
  \sum_{s\in\Edges}\int_s\nabla\Bu:\jumpt{\Bvh}\ds =
  \sum_{s\in\Edges}\int_s\big(\nabla\Bu-\Pis{k-1}(\nabla\Bu)\big):
  \big(\jumpt{\Bvh}-\Pis{0}\jumpt{\Bvh}\big)\ds,
\end{align*}
for the velocity.
Using the Cauchy-Schwartz inequality and then applying component-wise
the approximation estimates of
Theorem~\ref{thm:polynomialApproximation} to bound each of the
resulting terms, we obtain, cf.~\cite{Ayuso-Lipnikov-Manzini:2016}
or~\cite{Crouzeix-Raviart:1973},
\begin{align*}
  \abs{\sum_{s\in\Edges}\int_s\nabla\Bu:\jumpt{\Bvh}\ds\,} +
  \abs{\sum_{s\in\Edges}\int_s p\jump{\Bvh}\ds\,}
  \leq Ch^{r}
  \sum_{s\in\Edges}\Big(\norm{\Bu}_{r+1,E^{+}\cup E^{-}}+\norm{p}_{r,E^{+} \cup E^{-}}\Big)\abs{\Bvh}_{1,E^{+} \cup E^{-}},
\end{align*}
where for each side $s$ the symbols $E^{+}$ and $E^{-}$ denote the two
elements sharing that side.
As the number of edges/faces is assumed to be uniformly bounded, the
required result follows with a positive constant $C$ independent of
$h$ and $\Bu$.
\ENDPROOF

\begin{lem}[non-conformity error]
  \label{lem:non-conformity:error}
  Let $k,m\geq1$, $r=\min(k,\reg)$ be integer numbers.
  Let
  $\Bu\in\big[H^1_0(\Omega)\cap H^{m+1}(\Omega)\big]^d$ and
  $p\in L^2(\Omega)\slash{\R}\cap H^{m}(\Omega)$
  be the velocity and pressure solution of
  problem~\eqref{eq:Stokes:PDE:weak:a}-\eqref{eq:Stokes:PDE:weak:b},
  with source term $\Bf\in\big[L^2(\Omega)\big]^d$ and homogeneous
  boundary condition $\Bg=0$ on $\Gamma$.
  Then, for every $\Bvh\in\VVhkz$, it holds that
  \begin{align*}
    \abs{\as(\Bu,\Bvh) - (\Bf,\Bvh)}
    \leq Ch^{r}\Big(\norm{\Bu}_{r+1}+\norm{p}_{r}\Big)\,\norm{\Bvh}_{1}
    + \abs{\bs(\Bvh,p)}.
  \end{align*}
\end{lem}
\BEGINPROOF
Test~\eqref{eq:Stokes:PDE:a} against $\Bvh\in\VVhkz$.
For $\reg\geq1$ we have at least $\Bu\in\big[H^2(\Omega)\big]^2$ and
$p\in H^1(\Omega)$ and we can apply the Green's identity.
By rearranging the summation on the internal sides and using the
definition of the jump operators given in
section~\ref{subsec:mesh:jumps:polynomials}, we obtain:
\begin{align}
  \as(\Bu,\Bvh) - (\Bf,\Bvh)
  = \as(\Bu,\Bvh) - (-\Delta\Bu + \D p,\Bvh)
  = \bs(\Bvh,p)
  + \sum_{s\in\Edges}\int_s\left(
  \nabla\Bu:\jumpt{\Bvh}+p\jump{\Bvh}
  \right)\,\ds.
  \label{eq:non-conformity}
\end{align}
The assertion of the lemma follows by applying the result of
Lemma~\ref{lem:estimates:of:jumps}.
\ENDPROOF

\subsection{Error estimate for the velocity}
\label{subsect:error:estimates:velocity}
Let
\begin{align}
  \BW := \Big\{\,\Bw\in\big[H^1_0(\Omega)\big]^d\,|\,\Div\Bw=0\,\Big\}
  \quad\textrm{and}\quad
  \Whk := \Big\{\,\Bwh\in\VVhkz\,|\,\PiE{k-1}\circ\Div\Bwh=0\,\Big\}.
  \label{eq:W:Wh:def}
\end{align}
Using these definitions,
problem~\eqref{eq:Stokes:PDE:weak:a}-\eqref{eq:Stokes:PDE:weak:b} is
equivalent to~\cite{Crouzeix-Raviart:1973}:
\begin{align}
  \mbox{\emph{Find $\Bu\in\BW$ such that~~}}
  \as(\Bu,\Bv) = (\Bf,\Bv) \quad\forall\Bv\in\BW,
  \label{eq:Stokes:weak:velocity}
\end{align}
and
problem~\eqref{eq:Stokes:dPDE:weak:a}-\eqref{eq:Stokes:dPDE:weak:b} is
equivalent to~\cite{Crouzeix-Raviart:1973}:
\begin{align}
  \mbox{\emph{Find $\Buh\in\Whk$ such that~~}}
  \ah(\Buh,\Bvh) = (\Bfh,\Bvh) \quad\forall\Bvh\in\Whk.
  \label{eq:Stokes:VEM:velocity}
\end{align}

\medskip
\noindent
\begin{thm}[$H^1$ abstract a priori error bound for the velocity]
  \label{thm:vemStrang}
  Let $\Bu\in\BW$ be the solution of
  problem~\eqref{eq:Stokes:weak:velocity}, and $\Buh\in\Whk$ the
  solution of problem~\eqref{eq:Stokes:VEM:velocity} with $k\geq 1$.
  Then, it holds that:
  \begin{align}
    \bSnorm{\Bu-\Buh}
    \leq
    & \frac{1}{\alpha_*}\Bigg[\,
    (1+\alpha^*)
    \inf_{\Bvh\in\Whk}\bSnorm{\Bu-\Bvh} +
    (1+\alpha^*)
    \inf_{\Bq\in[\Phi^k_h]^d}\bSnorm{\Bu-\Bq}\nonumber\\[0.5em]
    &\quad\qquad+
    \sup_{\substack{\Bwh\in\Whk\\ \Bwh\neq 0}}\frac{\abs{ (\Bfh,\Bwh) - (\Bf,\Bwh) }}{\bSnorm{\Bwh}}
     \,+\sup_{\substack{\Bwh\in\Whk\\ \Bwh\neq 0}}\frac{\abs{\as(\Bu,\Bwh)-(\Bf,\Bwh)}}{\bSnorm{\Bwh}}
    \,\Bigg].
    \label{eq:vemStrang}
  \end{align}
  The last term in the right-hand side of the above error estimate
  measures the non-conformity error, i.e. it is non-zero because
  $\Whk$ is a non-conforming space.
\end{thm}
\BEGINPROOF
Let $\Bvh$ be an arbitrary element of $\Whk$ and let $\Bwh=\Buh-\Bvh$.
From the stability property of the virtual bilinear
form~\eqref{eq:stab:conds} and equation~\eqref{eq:Stokes:VEM:velocity}
it follows that
\begin{align*}
  \alpha_*\norm{\Buh - \Bvh}_{1,h}^2
  \leq \ah(\Buh-\Bvh,\Bwh)
  = (\Bfh,\Bwh) - \ah(\Bvh,\Bwh).
\end{align*}
Then, we add and subtract $(\Bf,\Bwh)$ and $\as(\Bu,\Bwh)$ and we
obtain:
\begin{align}
  \alpha_*\norm{\Buh - \Bvh}_{1,h}^2
  \leq
  \Big[ (\Bfh,\Bwh) - (\Bf,\Bwh) \Big] +
  \Big[ (\Bf, \Bwh) - \as(\Bu,\Bwh) \Big] +
  \Big[ \as(\Bu,\Bwh) - \ah(\Bvh,\Bwh) \Big].
  \label{eq:abstract:proof:velocity}
\end{align}
The first term on the right-hand side characterizes the approximation
of the source term $\Bf$ by $\Bfh$; the second term is the conformity
error determined by choosing the test function $\Bwh$ in $\Whk$
instead of $\BW$; the third term simultaneously characterizes the
approximation of $\Bu$ by $\Bvh$ and $\as(\cdot,\cdot)$ by
$\ah(\cdot,\cdot)$ in $\Whk$.
To separate this linked dependence, we reformulate the last term as
the summation of local contributions; then, we add and subtract
$\aE(\Bq,\Bwh)$ and $\ahE(\Bq,\Bwh)$ to each summation argument, where
$\Bq$ is a generic vector-valued function in $\big[\Phi^k_h\big]^d$,
and we obtain:
\begin{align*}
  \as(\Bu,\Bwh) - \ah(\Bvh,\Bwh)
  &=\sum_{E\in\Th}\Big[ \aE(\Bu,\Bwh) - \ahE(\Bvh,\Bwh) \Big]\\
  &=\sum_{E\in\Th}\Big[ \aE(\Bu-\Bq,\Bwh) +\big(\aE(\Bq,\Bwh) - \ahE(\Bq,\Bwh)\big) + \ahE(\Bq-\Bvh,\Bwh) \Big].
\end{align*}
Since $\Bq_{|E}$ is a polynomial vector of degree $k$, the
intermediate term above is zero due to the polynomial consistency
relation~\eqref{eq:pcons:conds}.
Furthermore, we transform the last term above by using the continuity
and stability of the virtual bilinear form:
\begin{align*}
  \ahE(\Bq-\Bvh,\Bwh)
  &
  \leq
  \left(\ahE(\Bq-\Bvh,\Bq-\Bvh)\right)^{\frac{1}{2}}
  \left(\ahE(\Bwh,\Bwh)\right)^{\frac{1}{2}}
  \leq \alpha^*
  \left(\aE(\Bq-\Bvh,\Bq-\Bvh)\right)^{\frac{1}{2}}
  \left(\aE(\Bwh,\Bwh)\right)^{\frac{1}{2}}
  \\[0.5em]
  &=
  \alpha^*\bSnormE{\Bq-\Bvh}\bSnormE{\Bwh}.
\end{align*}
Hence, for all $\Bvh, \Bwh\in\Whk$ with $\Bwh\neq\mathbf{0}$ and
$\Bq\in\big[\Phi_{h}^{k}\big]^d$, using the last inequality and the
continuity of $\aE(\cdot,\cdot)$, multiplying and dividing by
$\bSnorm{\Bwh}$, and using the triangle inequality
$\bSnorm{\Bq-\Bvh}\leq\bSnorm{\Bq-\Bu}+\bSnorm{\Bu-\Bvh}$, we find
that
\begin{align*}
  \alpha_*\bSnorm{\Buh-\Bvh}
  \leq
  \frac{\abs{(\Bfh,\Bwh) - (\Bf,\Bwh)}}{\bSnorm{\Bwh}} +
  \frac{\abs{(\Bf,\Bwh)-\as(\Bu,\Bwh)}}{\bSnorm{\Bwh}}
  + \alpha^*\norm{\Bu-\Bvh}_{1,h}+ \big(1+\alpha^*\big)\bSnorm{\Bu-\Bq}.
\end{align*}
The result now follows by applying the estimates above to the triangle
inequality
\begin{align*}
  \bSnorm{\Bu-\Buh}\leq\bSnorm{\Bu-\Bvh}+\bSnorm{\Bvh-\Buh},
\end{align*}
properly taking the supremum upper bound of the terms with $\Bwh$ and
the infimum over the arbitrary vector functions $\Bq$ and $\Bvh$.
\ENDPROOF

\medskip
\noindent
\begin{thm}[$H^1$ error bound for the velocity]
  \label{thm:h1ErrorBound}
  Let $k,\reg\geq1$, $r=\min(k,\reg)$ be integer numbers.
  Assuming that $\Bf\in\big[H^{\reg-1}(\Omega)\big]^d$, let
  $\Bu\in\big[H^1_0(\Omega)\cap H^{\reg+1}(\Omega)\big]^d$ and $p\in
  L^2(\Omega)\slash{\R}\cap H^{\reg}(\Omega)$ be the exact velocity
  and pressure solution to
  problem~\eqref{eq:Stokes:PDE:weak:a}-\eqref{eq:Stokes:PDE:weak:b}
  (with $\Bg=0$ on $\Gamma$).
  Let $(\Bfh,\Bvh):=\sum_{E\in\Th}(\Bfh,\Bvh)_E$, with ${\Bfh}_{|E}$
  defined as in section~\ref{subsec:rhs:approximation}.
  Denote by $\Buh\in\Whk$ the virtual element solution to
  problem~\eqref{eq:Stokes:VEM:velocity} under
  Assumption~\ref{assum:mesh:regularity} (mesh regularity), where
  $\Whk$ is the divergence-free non-conforming virtual element space
  of vector-valued functions defined in~\eqref{eq:W:Wh:def}.
  Then, there exists a constant $C$ independent of $h$ such
  that
  \begin{equation*}
    \bSnorm{\Bu-\Buh} \leq C h^r
    (\norm{\Bu}_{r+1} + \norm{p}_{r} + \norm{\Bf}_{r-1}).
  \end{equation*}
\end{thm}
\BEGINPROOF
Since $\Bu$ is also solution of
problem~\eqref{eq:Stokes:weak:velocity}, we start the proof of the
theorem by separately bounding the terms of the abstract bound of
Theorem~\ref{thm:vemStrang}.
The first term on the right-hand side of~\eqref{eq:vemStrang},
i.e. $\inf_{\Bvh\in\Whk}\bSnorm{\Bu-\Bvh}$, is easily bounded by
introducing the virtual interpolant $\BuI\in\VVhk$ of $\Bu$ as in
Theorem~\ref{thm:spaceApproximation} and noting that $\BuI$ belongs to
$\Whk$.
Indeed, for any element $E$, we have $\PiE{k-1}(\Div\BuI)\in\PE{k-1}$
and for any polynomial $q\in\PE{k-1}$ we have that:
\begin{align*}
  \int_{E}q\PiE{k-1}(\Div\BuI)\,\dx
  &
  = \int_{E}q\Div\BuI \dx
  = -\int_{E}\nabla q\cdot\BuI\dx + \sum_{s\in\d E}\int_sq\Bn_s\cdot\BuI\ds \\[0.5em]
  &= -\int_{E}\nabla q\cdot\Bu \dx + \sum_{s\in\d E}\int_sq\Bn_s\cdot\Bu \ds
  = \int_{E}q\Div\Bu\, \dx = 0.
\end{align*}
Hence, $\PiE{k-1}(\Div\BuI)=0$.
Likewise, the second term is bounded by using in each cell $E$ the
results of Theorem~\ref{thm:polynomialApproximation} since
\begin{align*}
  \inf_{\Bq\in[\Phi^k_h]^d}\bSnorm{\Bu-\Bq}\leq\bSnorm{\Bu-\PiE{k}(\Bu)}.
\end{align*}
The bound on the third term is given by Lemma~\ref{lem:rhs:error}
since $\Bwh\in\Whk\subset\VVhk$.
The bound on the fourth term is given by
Lemma~\ref{lem:non-conformity:error} and noting that $\bs(\Bwh,p)=0$
for every $\Bwh\in\Whk$.
Finally, the assertion of the theorem follows by combining the bounds
derived above and noting that the inequality constant $C$ may depend
only on the stability constants $\alpha_*$ and $\alpha^*$, and the
mesh regularity constant $\rho$.
\ENDPROOF

\vspace{-\baselineskip}
\begin{remark}
  Optimal order error estimates for the velocity approximation in the
  $L^2$-norm can be derived by duality arguments, see
  always~\cite{Crouzeix-Raviart:1973}. However, a more accurate
  approximation of the forcing terms than that provided
  by~\eqref{eq:Bfh:def} would be needed for the case $k=1,2$. This can
  be obtained by following the approach
  in~\cite{BeiraodaVeiga-Brezzi-Marini:2013}.
\end{remark}

\vspace{0.75\baselineskip}

\subsection{Error estimate for the pressure}
\label{subsect:error:estimates:pressure}

Let $\pI$ denote the piecewise polynomial function that is defined on
each elements $E$ of mesh $\Th$ by the orthogonal projection of $p$ on
the space of polynomials of degree $k-1$; formally,
$(\pI)_{|E}=\PiE{k-1}(p)\in\PE{k-1}$.
The accuracy of this approximation is characterized by
Theorem~\ref{thm:polynomialApproximation}.
To ease the notation we will also use the symbol $p^{\INTP}_{E}$ to
denote the restriction $(\pI)_{|E}$.

\medskip
\begin{thm}[Abstract $L^2$ a priori error bound for the pressure]
  \label{thm:vemStrangL2}
  Let $k\geq1$ be an integer number.
  Let
  $\Bu\in\big[H^1_0(\Omega)\cap H^2(\Omega)\big]^d$
  and
  $p\in L^2(\Omega)\slash{\R}\cap H^1(\Omega)$
  be the exact velocity and pressure solution to
  problem~\eqref{eq:Stokes:PDE:weak:a}-\eqref{eq:Stokes:PDE:weak:b}
  with the homogeneous boundary condition $\Bg=0$ on $\Gamma$.
  Denote the projection of $p$ in $\Phk$ by $\pI$.
  Let $(\Bfh,\Bvh):=\sum_{E\in\Th}(\Bfh,\Bvh)_E$, with ${\Bfh}_{|E}$
  defined as in section~\ref{subsec:rhs:approximation}.
  Let $\Buh\in\VVhkz$ and $\ph\in\Phk$ be the virtual element velocity
  and pressure solution to
  problem~\eqref{eq:Stokes:dPDE:weak:a}-\eqref{eq:Stokes:dPDE:weak:b}
  under Assumption~\ref{assum:mesh:regularity}.
  Then, it holds that
  \begin{align}
    &\norm{\ph-\pI}_0
    \leq
    \frac{\alpha^*}{\sqrt{\alpha_*}\beta}
    \left(\,
      \alpha^*
      \bSnorm{\Bu-\Buh} +
      (1+\alpha^*)
      \inf_{\Bq\in[\Phi^k_h]^d}\bSnorm{\Bu-\Bq}
      \phantom{
        \sup_{\substack{\Bvh\in\VVhkz\\ \Bvh\neq 0}}
        \frac{\Abs {\sum_{s\in\Edges}\int_sp\jump{\Bvh}\ds}}{\bSnorm{\Bvh}}
      }
    \right.
    \nonumber\\
    &\quad+
    \left.
    \sup_{\substack{\Bvh\in\VVhkz\\ \Bvh\neq 0}}
    \frac{\displaystyle
      \abss{ (\Bfh,\Bvh) - (\Bf,\Bvh) } +
      \abss{ \bs(\Bvh,p-\pI) } +
      \Abs {\sum_{s\in\Edges}\int_s\nabla\Bu:\jumpt{\Bvh}\ds} +
      \Abs {\sum_{s\in\Edges}\int_sp\jump{\Bvh}\ds}
    }{\bSnorm{\Bvh}}
    \right).
  \end{align}
\end{thm}
\BEGINPROOF
Take $q=\ph-\pI\in\Phk$.
From Lemma~\ref{lem:inf-sup}, we know that there exists a vector
$\Bvh$ such that
$\beta\norm{\ph-\pI}_{0}\,\TNORM{\Bvh}\leq\bh(\Bvh,\ph-\pI)$, and, in
view of the norm equivalence~\eqref{eq:equiv:norms}, it holds that
\begin{align}
  \beta\norm{\ph-\pI}_{0}
  \leq\frac{\bh(\Bvh,\ph-\pI)}{\sqrt{\alpha_*}\bSnorm{\Bvh}}
  \leq\sup_{\substack{\Bvh\in\VVhkz\\ \Bvh\neq 0}}
  \frac{\bh(\Bvh,\ph-\pI)}{\sqrt{\alpha_*}\bSnorm{\Bvh}}.
  \label{eq:pressure:proof:begin}
\end{align}
We use~\eqref{eq:Stokes:dPDE:weak:a}, and add and subtract
$(\Bf,\Bvh)$ and $\as(\Bu,\Bvh)$ to obtain:
\begin{align}
  \bh(\Bvh,\ph-\pI)
  &= (\Bfh,\Bvh) - \ah(\Buh,\Bvh) - \bh(\Bvh,\pI)
  \nonumber\\[0.5em]
  &=
  \Big[ (\Bfh,\Bvh) - (\Bf,\Bvh) \Big] +
  \Big[ (\Bf,\Bvh) - \as(\Bu,\Bvh) \Big] +
  \Big[ \as(\Bu,\Bvh) - \ah(\Buh,\Bvh) \Big]
  - \bh(\Bvh,\pI).
  \label{eq:pres:err:05}
\end{align}
We test~\eqref{eq:Stokes:PDE:a} against $\Bvh\in\VVhkz$, apply the
Green's identity, rearrange the summation on edges/faces in
$\mathcal{E}$, and introduce the jump notation
(cf. section~\ref{subsec:mesh:jumps:polynomials}) to obtain:
\begin{align}
  &(\Bf,\Bvh) - \as(\Bu,\Bvh)
  = (-\Delta \Bu + \D p,\Bvh) - \as(\Bu,\Bvh)
  \nonumber\\[0.5em]
  &\qquad = \bs(\Bvh,p)
  - \left(
    \sum_{E\in\Th}\int_{\d E}(\Bn\cdot\D)\Bu\cdot\Bvh\ds -
    \sum_{E\in\Th}\int_{\d E} p\Bn\cdot\Bvh\ds
  \right)
  \nonumber\\[0.5em]
  &\qquad = \bs(\Bvh,p)
  - \sum_{s\in\Edges} \int_s\Big((\Bn_s^+\cdot\nabla)\Bu\cdot\Bvh^++(\Bn_s^-\cdot\nabla)\Bu\cdot\Bvh^-\Big)\,\ds
  + \sum_{s\in\Edges} \int_s p(\Bn_s^+\cdot\Bvh^++\Bn_s^-\cdot\Bvh^-)\,\ds
  \nonumber\\[0.5em]
  &\qquad = \bs(\Bvh,p)
  - \sum_{s\in\Edges}\int_s\left(
    \nabla\Bu:\jumpt{\Bvh} - p\jump{\Bvh}
  \right)\,\ds.
  \label{eq:pres:err:10}
\end{align}
Also, from Remark~\ref{rem:bhE:bE} we know that
$\bh(\Bvh,\pI)=\bs(\Bvh,\pI)$, and using this relation and
\eqref{eq:pres:err:10} in~\eqref{eq:pres:err:05} we obtain:
\begin{align}
  \bh(\Bvh,\ph-\pI)
  &
  = \Big[ (\Bfh,\Bvh) - (\Bf,\Bvh) \Big]
  + \Big[ \as(\Bu,\Bvh) - \ah(\Buh,\Bvh) \Big]
  + \bs(\Bvh, p-\pI)
  \nonumber\\[0.5em]
  &
  -\sum_{s\in\Edges}\int_s\left(
    \nabla\Bu:\jumpt{\Bvh}-p\jump{\Bvh}
  \right)\,\ds.
  \label{eq:pres:err:20}
\end{align}
The second term on the right-hand side can be further bounded by
introducing a generic vector $\Bq\in[\Phi^k_h]^d$ and reasoning as in
the proof of Theorem~\ref{thm:vemStrang}, thus yielding
\begin{align}
  \abss{\as(\Bu,\Bvh) - \ah(\Buh,\Bvh)}\leq
  \left(\alpha^*\bSnorm{\Bu-\Buh}+ \big(1+\alpha^*\big)\bSnorm{\Bu-\Bq}\right)\bSnorm{\Bvh}.
  \label{eq:pres:err:25}
\end{align}
The assertion of the theorem follows
from~\eqref{eq:pressure:proof:begin} by taking the absolute value
of~\eqref{eq:pres:err:20} with~\eqref{eq:pres:err:25} $\bSnorm{\Bvh}$
and the supremum on $\Bvh$ and $\Bq$.
\ENDPROOF

\begin{thm}[$L^2$ a priori error bound for the pressure]
  Under the assumptions and notations of
  Theorem~\ref{thm:h1ErrorBound}, let $\Buh\in\VVhkz$ and $\ph\in\Phk$
  be the virtual element velocity and solution to
  problem~\eqref{eq:Stokes:dPDE:weak:a}-\eqref{eq:Stokes:dPDE:weak:b}.
  Then, there exists a constant $C>0$ depending only on the stability
  constants $\alpha_*$ and $\alpha^*$, the inf-sup constant $\beta$
  and the mesh regularity constant $\rho$ such that
  \begin{align}
    \norm{p-\ph} \leq C h^r
    \big(\norm{\Bu}_{r+1} + \norm{p}_{r} + \norm{\Bf}_{r-1}\big).
  \end{align}
  \vspace{-0.5cm}
\end{thm}
\BEGINPROOF
The proof is just a matter of bounding the terms of the abstract bound
of Theorem~\ref{thm:vemStrangL2}.
Bounds for the first, second, and third term are already given in
Theorem~\ref{thm:h1ErrorBound}.
The bound for the term containing $\bs(\Bvh,p-\pI)$ follows from the
polynomial approximation results of
Theorem~\ref{thm:polynomialApproximation}.
The last two terms are bounded by using
Lemma~\ref{lem:estimates:of:jumps}.
\ENDPROOF

\section{Implementation details}
\label{sec:implementation}

\newcommand{\matPk} {\bm\Pi_{k}^{E}}
\newcommand{\matPkk}{\bm\Pi_{k-1}^{E,\xl}}
\newcommand{\matPkkp}{\bm\Pi_{k-1}^{E,x_{l'}}}
\newcommand{\vphi}  {\phi}
\newcommand{\matPkN} {\bm\Pi_{k}^{\nabla}}

\newcommand{\xl} {x^l}
\newcommand{\uhl}{u_h^l}
\newcommand{\vhl}{v_h^l}
\newcommand{\uhlp}{u_h^{l'}}
\newcommand{\vhlp}{v_h^{l'}}

\newcommand{\fl}{f^l}
\newcommand{\fhl}{f_h^l}

\newcommand{\coeff}[2]{c^{#1}_{#2}}

According to section~\ref{sect:non-conforming:virtual:element:space},
the key component of the VEM implementation is the construction of the
projection operators $\PiE{k-1}\circ\nabla$, $\PiE{k-1}\circ\Div$, and
$\Pi^{\nabla}_{k}$.
Once these projection operators are constructed, all terms in the
local bilinear forms can be computed as integrals of polynomials just
as in the standard FEM with the only exception of the stabilization
term in $\ahE$.
This latter term does not require any integration but is directly
defined through the action of the projector operator
$\Pi^{\nabla}_{k}$ on the degrees of freedom.

The matrix representation of $\PiE{k-1}\circ\nabla$ and
$\PiE{k-1}\circ\Div$ can be derived from the corresponding representation 
of the scalar projection operator $\PiE{k-1}\circ\tfrac{\d}{\d x_i}$,
$i=1,\ldots,d$,
which has already been worked out in~\cite{Cangiani-Manzini-Sutton:2016}.
Similarly, we refer to~\cite{Ahmad-Alsaedi-Brezzi-Marini-Russo:2013} for the details on the computation of $\Pi^{\nabla}_{k}$.
In the rest of this section we show how to compute the two terms of
the bilinear form $\ahE$, i.e., the consistency and stability term,
cf. Remark~\ref{remark:consistency:stability}, assuming a matrix 
representation of the projectors is given.
The implementation formulas for the bilinear form $\bhE$ and the
right-hand side linear functional $(\Bfh,\cdot)_{E}$ can be 
derived similarly and are not shown here.

Consider the Lagrangian basis $\{\vphi_i\}$ for the scalar virtual
element space $\Vhk(E)$ associated with the degrees of freedom
introduced in
Section~\ref{sect:scalar:non-conforming:virtual:element:space}.
We collect the coefficients of the expansion of each monomials $\ma$
on the basis $\big\{\vphi_i\big\}$ for $\alpha=1,\ldots,N_{d,k}$ on
the columns of matrix $\mD$, so that
$\ma=\sum_{j=1}^{N_E}\vphi_j\mD_{j,\alpha}$.
Similarly, we collect the coefficients of the expansions of the
polynomials $\Pi^{\nabla}_{k}\vphi_i$ and
$\PiE{k-1}\big({\d\vphi_i}\slash{\d x_i}\big)$ with respect to the
monomial basis $\big\{\ma\big\}$ on the columns of matrices $\matPkN$
and $\matPkk$, respectively.
Hence, these projections can be expressed by the formulas
\begin{align}
  \PiE{k-1}\left(\frac{\partial\vphi_i}{\partial x_l}\right) &= \sum_{\alpha=1}^{N_{d,k}}\ma \big(\matPkk\big)_{\alpha,i}
  \quad\textrm{and}\quad
  \Pi^{\nabla}_{k}\vphi_i &= \sum_{\alpha=1}^{N_{d,k}}\ma \big(\matPkN\big)_{\alpha,i} = \sum_{j=1}^{N_E}\vphi_j \big(\mD\matPkN\big)_{j,i},
  \label{eq:proj:dphi:dxl}
\end{align}
where $( \star )_{i,j}$ denotes the $i,j$-th element of a given matrix
argument $\star$.
Then, for the local discrete velocity space $\VVhk(E)$, we consider
the set of basis functions $\{\bm{\phi}_i^l\}$ that are such that
$(\bm{\phi}_i^l)_{l'}=\phi_i$ if $l'=l$ and $(\bm{\phi}_i^l)_{l'}=0$
otherwise.
The generic entry of the consistency term of $\ahE$ in~\eqref{eq:ahE} is given by
\begin{align}
  &\int_E\PiE{k-1}(\D\bm{\phi}_i^l):\PiE{k-1}(\D\bm{\phi}_j^l)\,\dx
  = \sum_{l'=1}^d\int_{E}\PiE{k-1}\left(\frac{\partial\phi_i}{\partial x_{l'}}\right)\,\PiE{k-1}\left(\frac{\partial\phi_j}{\partial x_{l'}}\right)\,\dx\nonumber\\[0.5em]
  &\quad= \sum_{l'=1}^{d}\sum_{\alpha,\beta=1}^{N_{d,k-1}}\big(\matPkkp\big)_{\alpha,i}\big(\matPkkp\big)_{\beta,j}\int_{E}\ma\,\mb\,\dx
  = \sum_{l'=1}^{d} \Big(\big(\matPkkp\big)^T\mH \big(\matPkkp\big)\big)_{i,j},
\end{align}
where $\mH$ is the matrix with coefficients
$\big(\mH\big)_{\alpha,\beta}=\int_{E}\ma\,\mb\,\dx$ for
$\alpha,\beta=1,\ldots,N_{d,k-1}$.
These coefficients may be computed exactly in special cases or by
applying a sufficiently accurate integration rule in general.
Similarly, the generic entry of the stabilisation term of
$\ahE$, which is defined in~\eqref{eq:stabterm}, is given by
\begin{align}
  \shE\left(\,(I-\Pi^{\nabla}_{k})\bm{\phi}_i^l, (I-\Pi^{\nabla}_{k})\bm{\phi_j^l}\,\right) 
  &= \sum_{r=1}^{N_E}\bm{\chi}_r\left((I-\Pi^{\nabla}_{k})\vphi_i\right)\cdot\bm{\chi}_r\left((I-\Pi^{\nabla}_{k})\vphi_j\right)\nonumber\\[0.5em]
  &= \Big(\big(\mI-\mD\matPkN\big)^T\big(\mI-\mD\matPkN\big)\Big)_{i,j}
\end{align}
since, trivially, $\bm{\chi}_r(I-\Pi^{\nabla}_{k})\vphi_i=\big(\mI-\mD\matPkN\big)_{i,r}$.

\section{Numerical Results}
\label{sect:numerical:experiments}

\newcommand{\MeshONE}  {$\mathcal{M}_1$}
\newcommand{\MeshTWO}  {$\mathcal{M}_2$}
\newcommand{\MeshTHREE}{$\mathcal{M}_3$}
\newcommand{\HAT}[1]{\widehat{#1}}
\newcommand{\TABROW}[0]{}

The numerical experiments that we present in this section are aimed at
confirming the \emph{a priori} analysis developed in the previous
sections.
In a preliminary stage, the consistency of non-conforming VEM,
i.e. the exactness for polynomial solutions, has been tested
numerically by solving the Stokes equation with boundary and source
data determined by $\Bu(x,y)=(y^m,x^m)$ and $p=x^m+y^m$ on different
set of polygonal meshes and for $m=1$ to $4$.
In all the cases, we measure an error whose magnitude is of the order
of the arithmetic precision, thus confirming this property.

To study the accuracy of the method we solve the problem with the
following solution on the domain $\Omega=]0,1[\times]0,1[$:
\begin{align}
  \Bu(x,y) = \left(
    \begin{array}{c}
      2\pi\,f(x)\sin(2\pi y)\\
      f'(x)\cos(2\pi y)
    \end{array}
  \right),\qquad
  p(x,y) = \sin(2\pi x)\sin(2\pi y),
  \label{eq:exact:solution}
\end{align}
with $f(x)=x^5\,e^{-x}$.
The forcing term and the Dirichlet boundary condition are set in
accordance with~\eqref{eq:exact:solution}.

\begin{figure}[t]
  \centering
  \begin{tabular}{ccc}
    \includegraphics[width=0.28\textwidth]{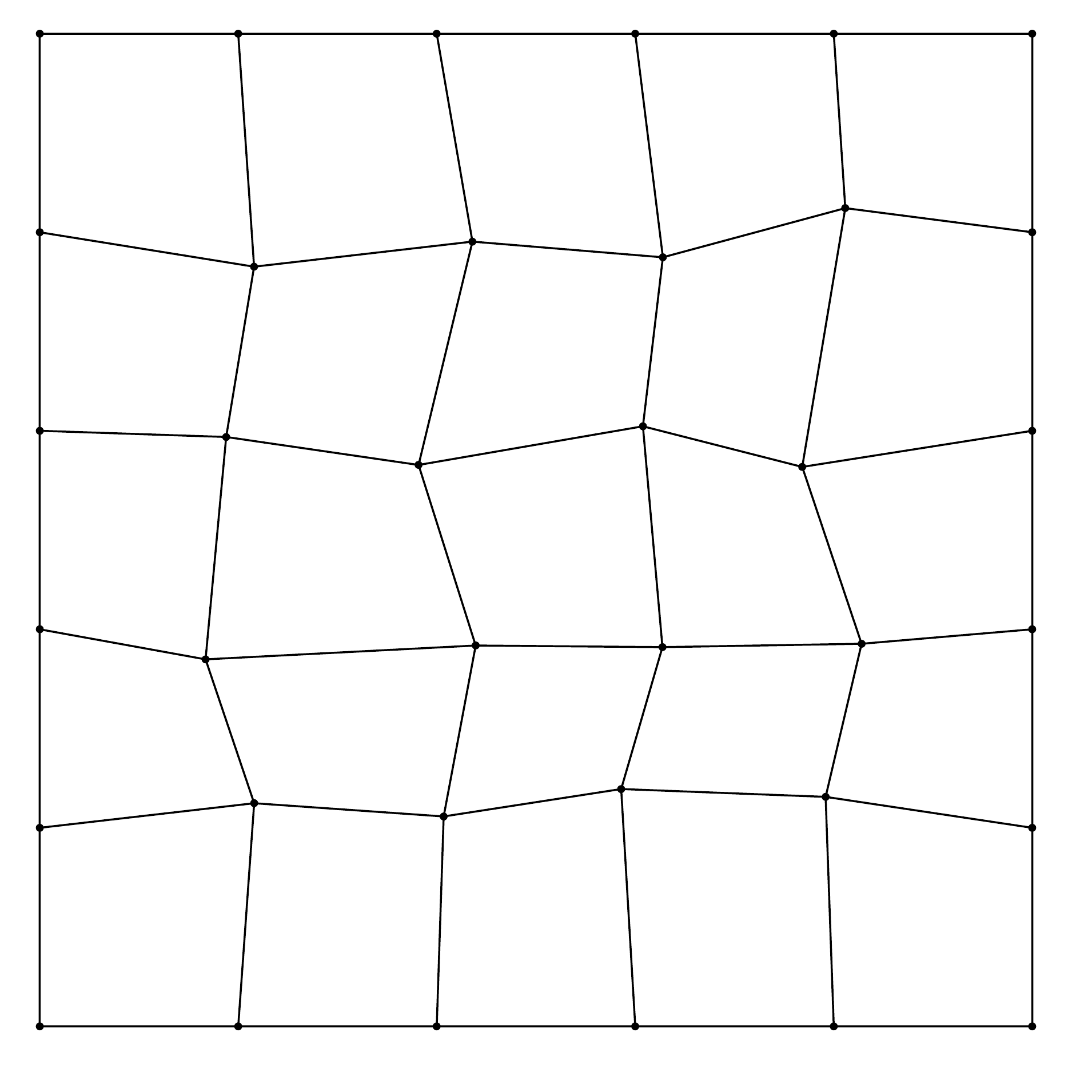}&\hspace{-0.25cm}
    \includegraphics[width=0.28\textwidth]{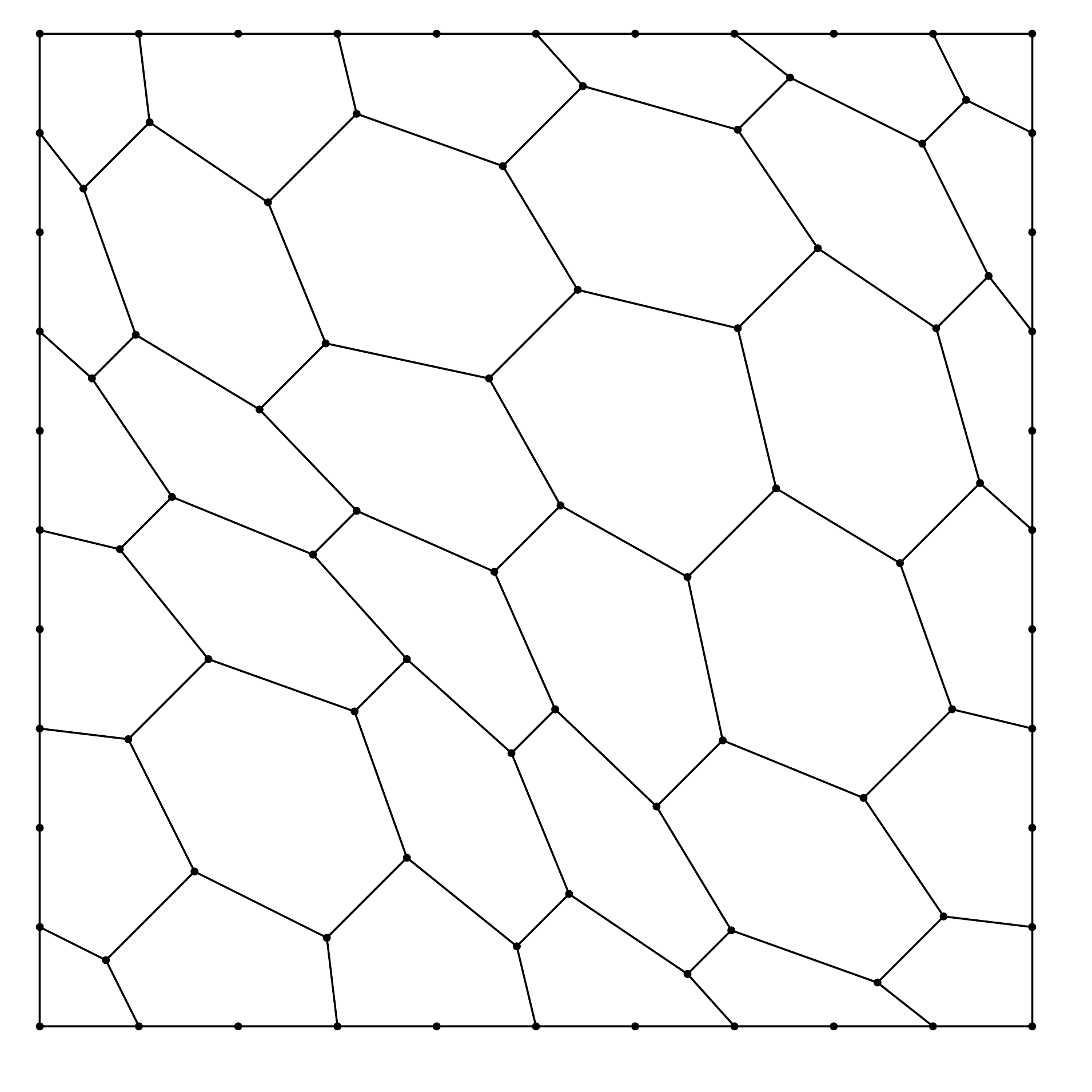}&\hspace{-0.25cm}
    \includegraphics[width=0.28\textwidth]{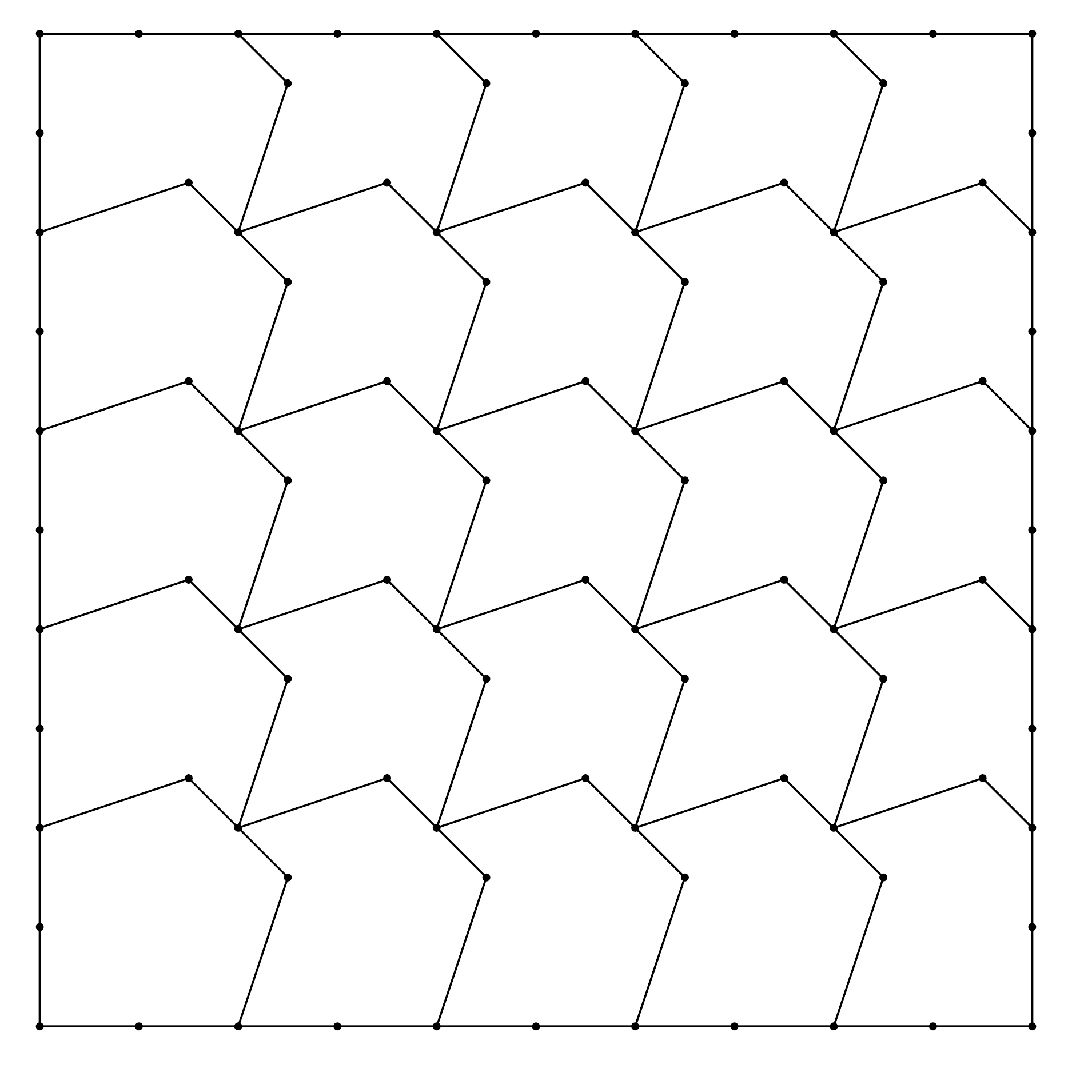}\\
    \includegraphics[width=0.28\textwidth]{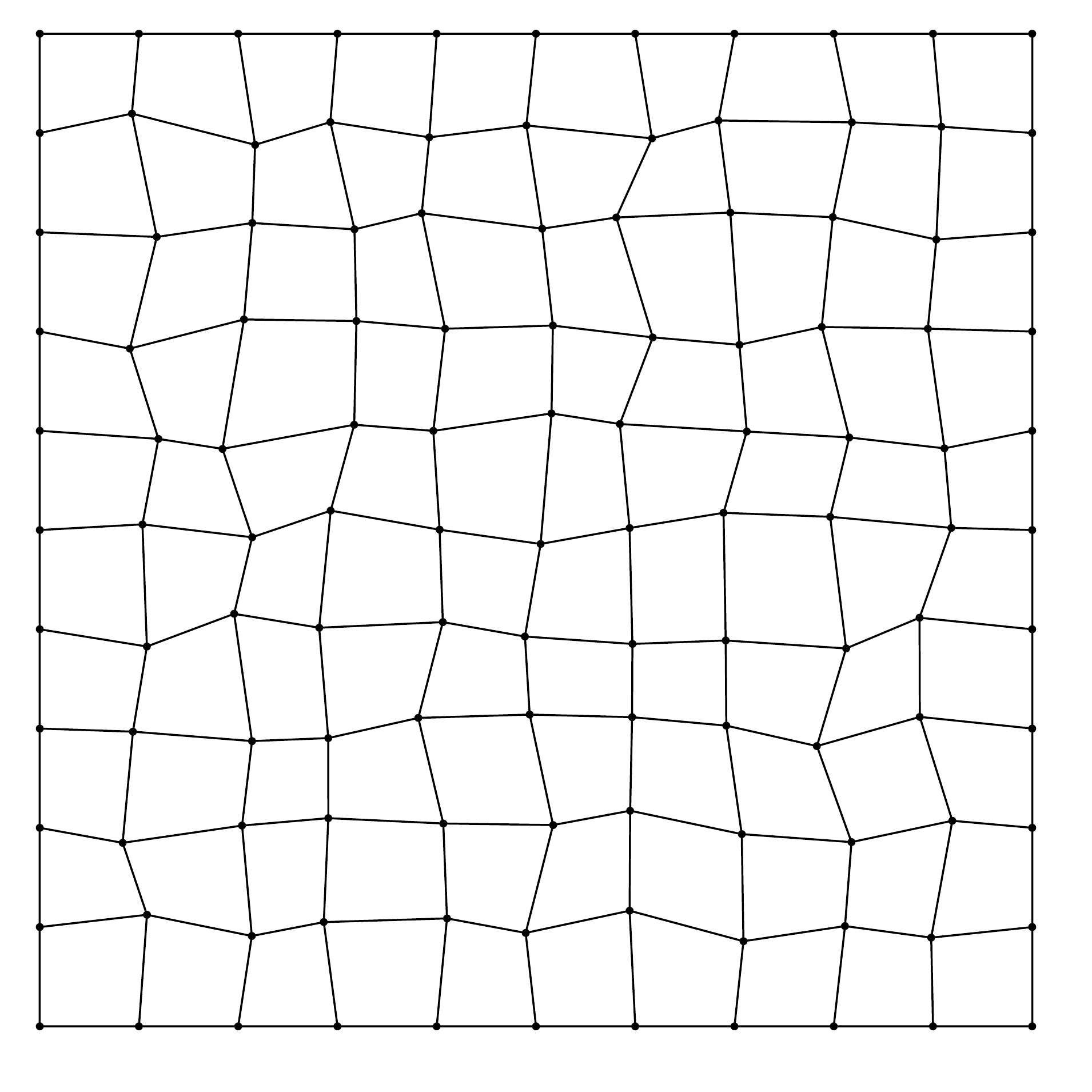}&\hspace{-0.25cm}
    \includegraphics[width=0.28\textwidth]{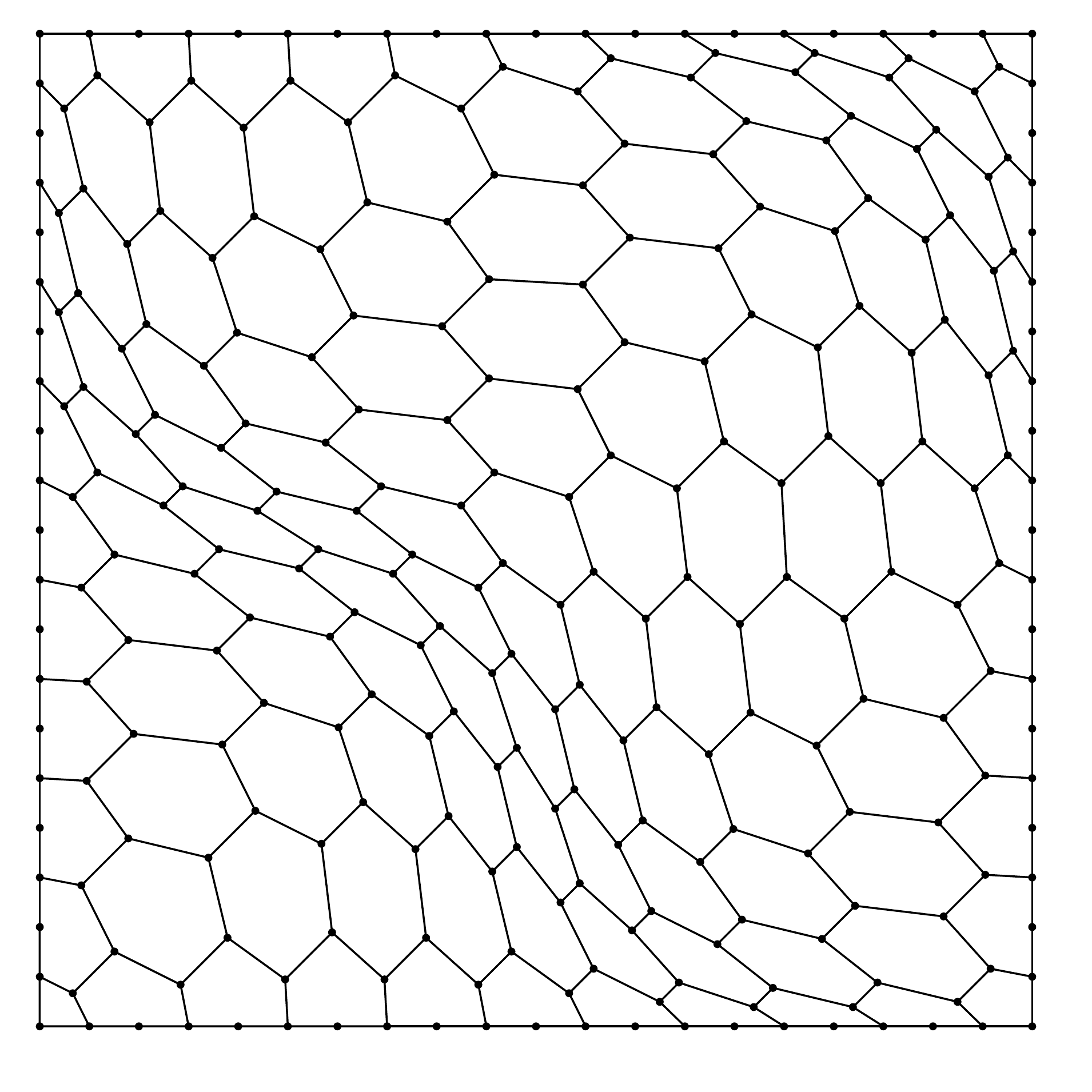}&\hspace{-0.25cm}
    \includegraphics[width=0.28\textwidth]{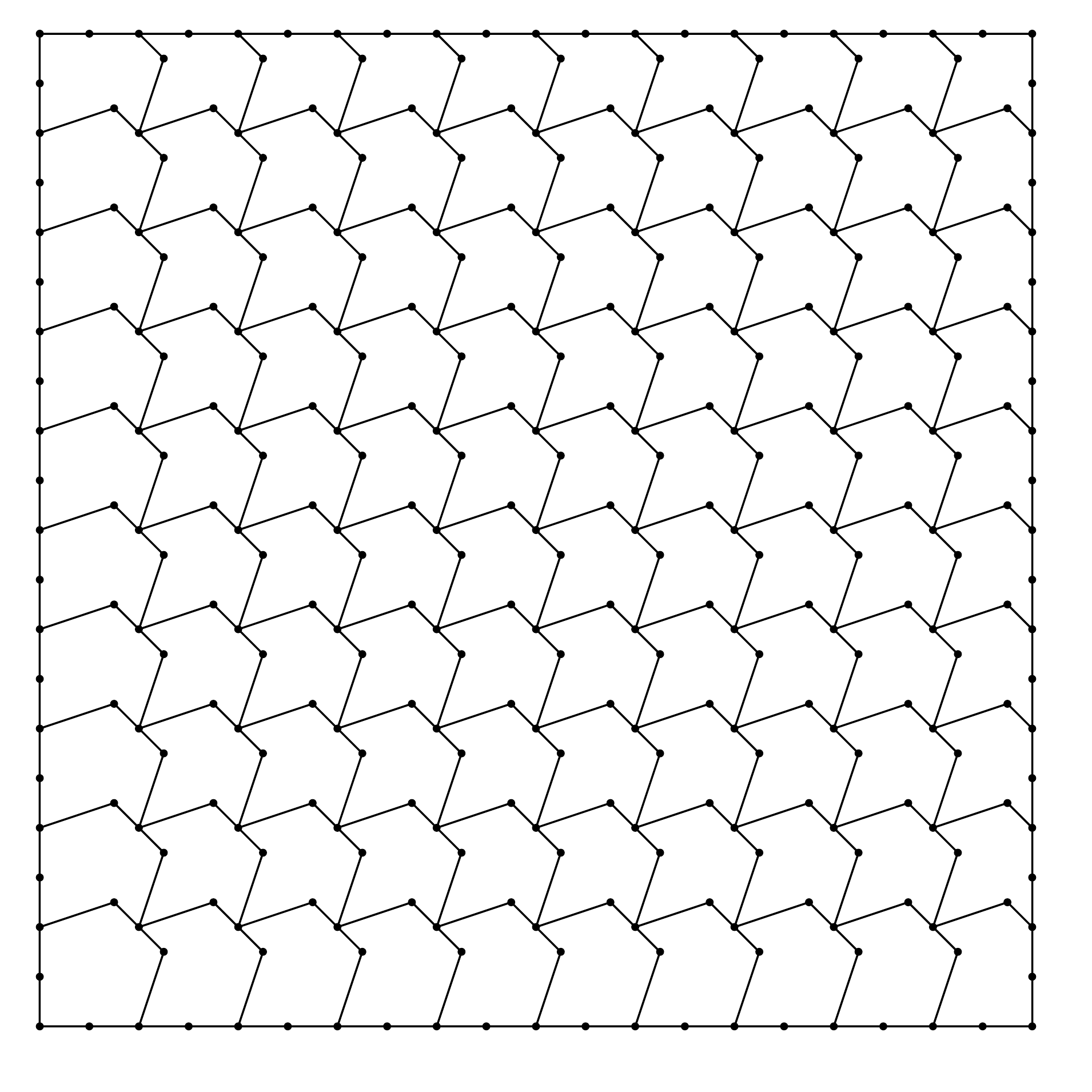}\\
  \end{tabular}
  \caption{First (top) and second (bottom) mesh of the three mesh families
  \MeshONE{} (left), \MeshTWO{} (middle) and \MeshTHREE{} (right).}
  \label{fig:meshes}
  \vspace{-0.125cm}
\end{figure}

The performance of the VEM for $k=1,2,3,4$ is investigated by
evaluating the rate of convergence on three different sequences of
five meshes, labeled by~\MeshONE{}, \MeshTWO{} and \MeshTHREE{},
respectively.
The top panels of Fig.~\ref{fig:meshes} show the first mesh of each
sequence and the bottom panels show the mesh of the first refinement.
The meshes in \MeshONE{} are built by partitioning the domain $\Omega$
into square cells and relocating each interior node to a random
position inside a square box centered at that node.
The sides of this square box are aligned with the coordinate axis and
their length is equal to $0.8$ times the minimum distance between two
adjacent nodes of the initial square mesh.
The meshes in \MeshTWO{} are built as follows.
First, we determine a primal mesh by remapping the position
$(\HAT{x},\HAT{y})$ of the nodes of a uniform square partition of
$\Omega$ by the smooth coordinate transformation:
\begin{align*}
  x &= \HAT{x} + (1\slash{10}) \sin(2\pi\HAT{x})\sin(2\pi\HAT{y}),\\
  y &= \HAT{y} + (1\slash{10}) \sin(2\pi\HAT{x})\sin(2\pi\HAT{y}).
\end{align*}
The corresponding mesh of \MeshTWO{} is built from the primal mesh by
splitting each quadrilateral cell into two triangles and connecting
the barycenters of adjacent triangular cells by a straight segment.
The mesh construction is completed at the boundary by connecting the
barycenters of the triangular cells close to the boundary to the
midpoints of the boundary edges and these latters to the boundary
vertices of the primal mesh.
The meshes in \MeshTHREE{} are obtained by filling the unit square
with a suitably scaled non-convex octagonal reference cell.

All the meshes are parametrised by the number of partitions in each
direction.
The starting mesh of every sequence is built from a $5\times 5$
regular grid, and the refined meshes are obtained by doubling this
resolution.
Mesh data for each refinement level, i.e., numbers of mesh elements,
number of edges, number of vertices, are reported in
Table~\ref{tab:mesh:data}.
More details on these mesh constructions can be found
in~\cite{BeiraodaVeiga-Manzini:2014,BeiraodaVeiga-Lipnikov-Manzini:2011}.
The mesh data structures are created and managed using the C++ mesh
manager tool described in~\cite{bertolazzi2002algorithm}.


\renewcommand{\TABROW}[5]{$#1$ & $#2$ & $#3$ & $#4$ & $#5$\\}

\newcommand{\ilev}{n}
\newcommand{\nE}{\mathcal{N}_{E}}
\newcommand{\nF}{\mathcal{N}_{e}}
\newcommand{\nv}{\mathcal{N}_{v}}

\begin{table}
  \begin{center}
    \begin{tabular}{|c|ccc|c|ccc|c|ccc|c|}
      \hline 
      &\multicolumn{4}{|c|}{Randomised quadrilaterals}
      &\multicolumn{4}{|c|}{Remapped hexagons}
      &\multicolumn{4}{|c|}{Non-convex octagons}\\
      \hline
      {$\ilev$}&{$\nE$}&{$\nF$}&{$\nv$}&{$h$}&{$\nE$}&{$\nF$}&{$\nv$}&{$h$}&{$\nE$}&{$\nF$}&{$\nv$}&{$h$}\\
      \hline
      {$1$}&{\!$  25$}&{\!$   60$}&{\!$  36$} &{$0.331$}
      {}   &{\!$  36$}&{\!$  125$}&{\!$   90$}&{$0.328$}
      {}   &{\!$  25$}&{\!$  120$}&{\!$   96$}&{$0.291$}\\
      {$2$}&{\!$ 100$}&{\!$  220$}&{\!$ 121$} &{$0.186$}
      {}   &{\!$ 121$}&{\!$  400$}&{\!$  280$}&{$0.185$}
      {}   &{\!$ 100$}&{\!$  440$}&{\!$  341$}&{$0.146$}\\
      {$3$}&{\!$ 400$}&{\!$  840$}&{\!$ 441$} &{$0.094$}
      {}   &{\!$ 441$}&{\!$ 1400$}&{\!$  960$}&{$0.097$}
      {}   &{\!$ 400$}&{\!$ 1680$}&{\!$ 1281$}&{$0.073$}\\
      {$4$}&{\!$1600$}&{\!$ 3280$}&{\!$1681$} &{$0.047$}
      {}   &{\!$1681$}&{\!$ 5200$}&{\!$ 3520$}&{$0.049$}
      {}   &{\!$1600$}&{\!$ 6560$}&{\!$ 4961$}&{$0.036$}\\
      {$5$}&{\!$6400$}&{\!$12960$}&{\!$6561$} &{$0.024$}
      {}   &{\!$6561$}&{\!$20000$}&{\!$13440$}&{$0.025$}
      {}   &{\!$6400$}&{\!$25920$}&{\!$19521$}&{$0.018$}\\
      \hline
    \end{tabular}
  \end{center}
  \caption{Mesh data for the meshes in \MeshONE{}, \MeshTWO{}, and \MeshTHREE{};
    $\nE$, $\nF$ and $\nv$ are the numbers of mesh elements, interfaces
    and vertices, respectively, and $h$ is the mesh size parameter.}
  \label{tab:mesh:data}
  \vspace{-0.25cm}
\end{table}

For the approximation of the velocity, we compare the polynomial
quantities $\PiE{k}(\Buh)$ and $\PiE{k-1}(\nabla\Buh)$ with the exact
velocity $\Bu$ and the gradient $\nabla\Bu$.
We recall that we can compute these projections \emph{exactly} using
only the degrees of freedom of the vector and scalar field $\Buh$
although we do not know this field.
For the approximation of the pressure, we compare the piecewise
polynomial fields $\ph$ and $\nabla\ph$ with the exact pressure $p$
and gradient $\nabla p$.

\begin{figure}[t]
  \hspace{0.5mm}
  \begin{tabular}{cc}
      \begin{overpic}[scale=0.73]{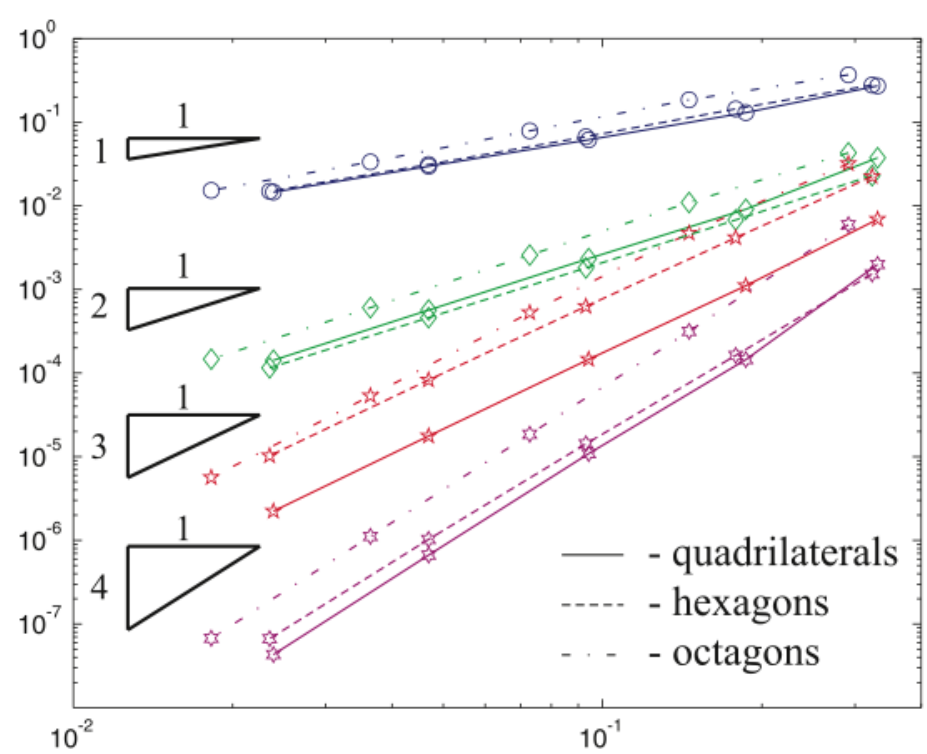}
        \put(75,-8){\textbf{Mesh size $h$}}
        \put(-12,38){\begin{sideways}\textbf{Pressure $L^2$ errors}\end{sideways}}
      \end{overpic}
      &\qquad
      \begin{overpic}[scale=0.73]{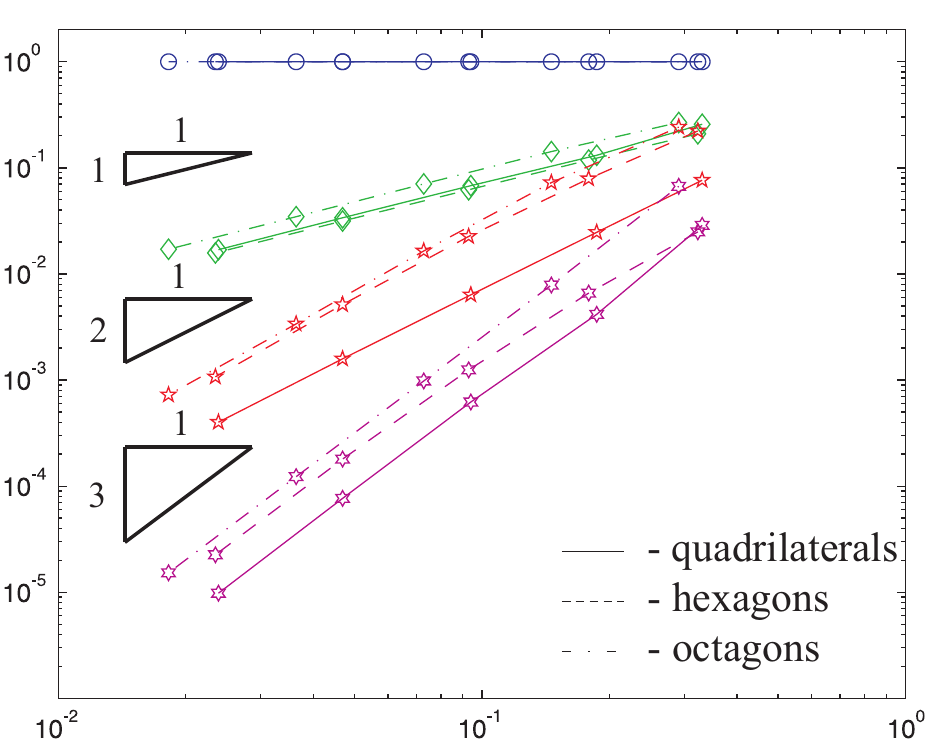}
        \put(75,-8){\textbf{Mesh size $h$}}
        \put(-12,38){\begin{sideways}\textbf{Pressure $H^1$ errors}\end{sideways}}
      \end{overpic}
  \end{tabular}
  \caption{
    Relative error curves for the non-conforming VEM approximation of the pressure (left)
    and its gradient (right) on the three mesh families \MeshONE{}-\MeshTHREE{}
    with $k$ from $1$ to $4$.
    The expected slopes are indicated by triangles and labels on the
    plots.  }
  \label{fig:errors:joint-meshes:pressure}
  \vspace{-0.25cm}
\end{figure}

\begin{figure}[t]
  \hspace{0.5mm}
  \begin{tabular}{cc}
      \begin{overpic}[scale=0.73]{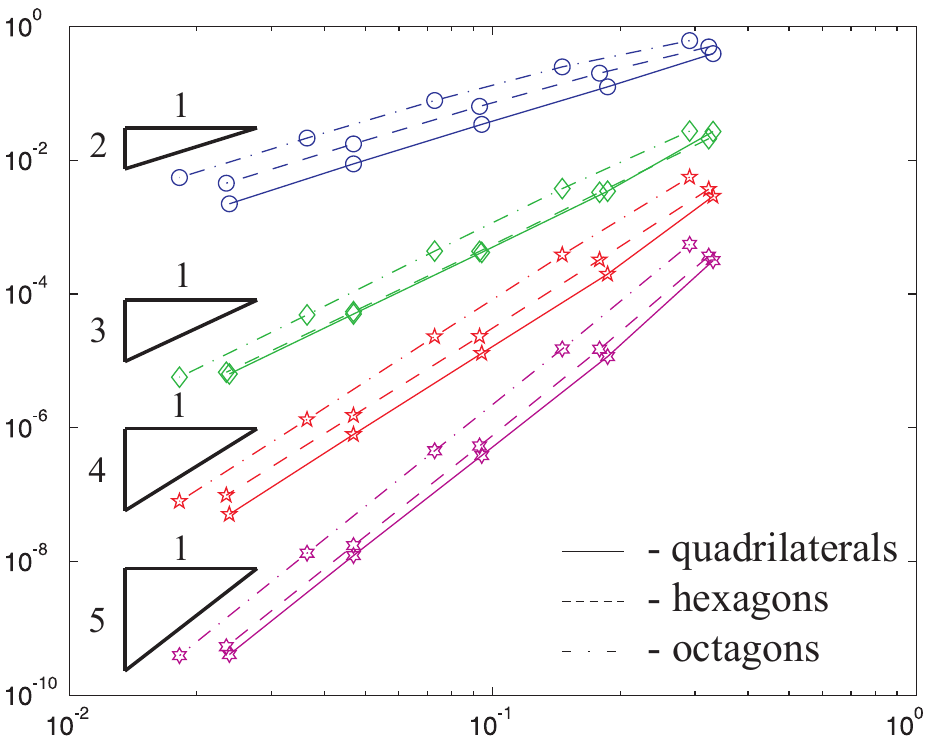}
        \put(75,-8){\textbf{Mesh size $h$}}
        \put(-12,38){\begin{sideways}\textbf{Velocity $L^2$ errors}\end{sideways}}
      \end{overpic}
      &\qquad
      \begin{overpic}[scale=0.73]{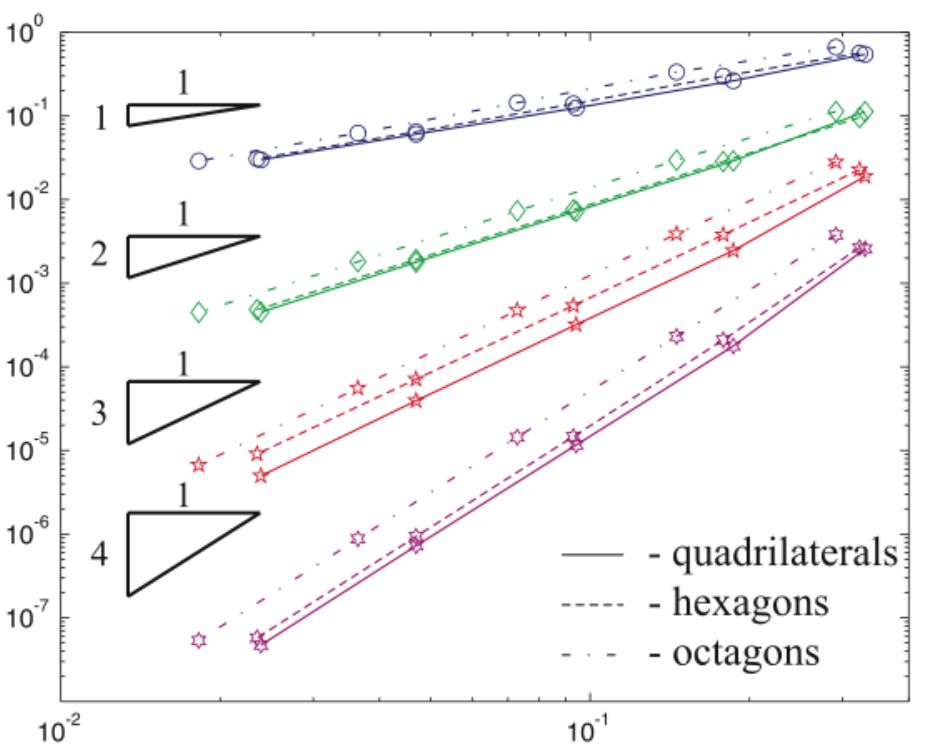}
        \put(75,-8){\textbf{Mesh size $h$}}
        \put(-12,38){\begin{sideways}\textbf{Velocity $H^1$ errors}\end{sideways}}
      \end{overpic}
  \end{tabular}
  \caption{
    Relative error curves for the non-conforming VEM approximation of the velocity (left)
    and its gradient (right) on the three mesh families \MeshONE{}-\MeshTHREE{}
    with $k$ from $1$ to $4$.
    The expected slopes are indicated by triangles and labels on the
    plots.  }
  \label{fig:errors:joint-meshes:velocity}
  \vspace{-0.25cm}
\end{figure}

The relative error curves for pressure and velocity versus the mesh
size $h$ are shown respectively in the log-log plots of
Figures~\ref{fig:errors:joint-meshes:pressure}
and~\ref{fig:errors:joint-meshes:velocity} for the three mesh
sequences as indicated therein.
The plots on the left show the relative errors for the approximation
of the velocity or pressure field, while the plots on the right show
the relative errors for the approximation of the field's gradient.
The expected slopes are shown for each error curve directly on the
plots and indicated by numerical labels.

These results are in very good agreement with the convergence rates
that are predicted by the analysis of the previous sections.

\section{Conclusions}
\label{sect:conclusion}
We presented the non-conforming formulation of the virtual element
method for the steady Stokes problem.
We have been able to construct approximations of any order in two and
three space dimensions in a unified fashion, a feat which is still out
of reach for standard non-conforming finite elements.
Moreover, the method is naturally defined on general polygonal and
polyhedral meshes.
In particular non-convex polygons and polyhedra with parallel adjacent
interfaces are allowed.
The formulation of the method relies on the element-wise construction
of virtual approximation spaces for velocity and pressure, which are
characterized by a positive integer order $k$.
The local approximation space for the velocity contains vectors of
polynomials of order $k$ plus other functions that are not computed
explicitly and dealt only in terms of their degrees of freedom.
The local approximation space for the pressure consists of all
polynomials of order $k-1$.
We proved that the velocity-pressure pair of global approximation spaces
satisfies the inf-sup condition, from which the stability
and well-posedness of the scheme follow.
We also proved the optimal convergence of the numerical approximation
to the velocity and pressure solution fields.
The accuracy of the approximation is determined by the degree $k$ of
the polynomials and optimal a priori error estimates was derived for
the velocity and pressure.


\vskip3mm
\section*{Acknowledgements}

The first author was partially supported by the Engineering and 
Physical Sciences Research Council of the United Kingdom (Grant
EP/L022745/1).
The second and third authors were partially supported by the
Laboratory Directed Research and Development program (LDRD),
U.S. Department of Energy Office of Science, Office of Fusion Energy
Sciences, under the auspices of the National Nuclear Security
Administration of the U.S. Department of Energy by Los Alamos National
Laboratory, operated by Los Alamos National Security LLC under
contract DE-AC52-06NA25396.
These supports are gratefully acknowledged.


\vskip4mm



\end{document}